\pdfoutput=1
\RequirePackage{ifpdf}
\ifpdf % We are running pdfTeX in pdf mode
\documentclass[pdftex]{sigma}
\else
\documentclass{sigma}
\fi

\numberwithin{equation}{section}

\usepackage{bbm}

\newcommand{\C}{\mathbb C}

\newcommand{\Z}{\mathbb Z}
\newcommand{\N}{\mathbb N}

\newcommand{\+}{\!+\!}

\newcommand{\Gt}{\mathfrak{N}}

\newcommand{\ddoAW}{\mathbb{D}}
\newcommand{\moAW}{\mathbb{M}}

\DeclareMathAlphabet{\mathpzc}{OT1}{pzc}{m}{it}

\begin{document}

\allowdisplaybreaks

\renewcommand{\PaperNumber}{097}

\FirstPageHeading

\ShortArticleName{Construction of a Lax Pair for the $E_6^{(1)}$ $q$-Painlev\'e System}

\ArticleName{Construction of a Lax Pair\\ for the $\boldsymbol{E_6^{(1)}}$ $\boldsymbol{q}$-Painlev\'e System}

\Author{Nicholas S.~WITTE~$^\dag$ and Christopher M.~ORMEROD~$^\ddag$}
\AuthorNameForHeading{N.S.~Witte and C.M.~Ormerod}

\Address{$^\dag$~Department of Mathematics and Statistics,
University of Melbourne, Victoria 3010, Australia}
\EmailD{\href{nsw@ms.unimelb.edu.au}{nsw@ms.unimelb.edu.au}}
\URLaddressD{\url{http://www.ms.unimelb.edu.au/~nsw/}}

\Address{$^\ddag$~Department of Mathematics and Statistics,
La Trobe University, Bundoora VIC 3086, Australia}
\EmailD{\href{C.Ormerod@latrobe.edu.au}{C.Ormerod@latrobe.edu.au}}

\ArticleDates{Received September 05, 2012, in f\/inal form November 29, 2012; Published online December 11, 2012}

\Abstract{We construct a Lax pair for the $ E^{(1)}_6 $ $q$-Painlev\'e system from f\/irst principles by employing the general theory of
semi-classical orthogonal polynomial systems characterised by divided-dif\/ference operators on
discrete, quadratic lattices [arXiv:1204.2328]. Our study treats one special case of such lattices~-- the $q$-linear lattice~--
through a natural generalisation of the big $q$-Jacobi weight.
As a by-product of our construction we derive the coupled f\/irst-order $q$-dif\/ference equations for the
$ E^{(1)}_6 $ $q$-Painlev\'e system, thus verifying our identif\/ication. Finally we establish the
correspondences of our result with the Lax pairs given earlier and separately by Sakai and Yamada, through explicit
transformations.}

\Keywords{non-uniform lattices; divided-dif\/ference operators; orthogonal polynomials; semi-classical weights; isomonodromic deformations; Askey table}
\Classification{39A05; 42C05; 34M55; 34M56; 33C45; 37K35}

\section{Background and motivation}\label{Start}

Since the recent discoveries of $ q $-analogues of the Painlev\'e equations, see for example~\cite{JS_1996} and~\cite{RGTT_2001}
which are of relevance to the present study, and their
classif\/ication (of these and others) according to the theory of rational surfaces by Sakai~\cite{Sa_2001}
interest has grown in f\/inding Lax pairs for these systems. This problem also has the independent interest as
a search for discrete and $ q $-analogues to the isomonodromic systems of the continuous Painlev\'e
equations, and an appropriate analogue to the concept of monodromy. Such interest, in fact, goes back to the
period when the discrete analogues of the Painlev\'e equations were f\/irst discussed, as one can see in~\cite{PNGR_1992}.

In this work we illustrate a general method for constructing Lax pairs for all the systems
in the Sakai scheme, as given in the study~\cite{Wi_2010a}, with the particular case of the $ E^{(1)}_6 $ system.
In this method all aspects of
the Lax pairs are constructed, and in the end we verify the identif\/ication with the $ E^{(1)}_6 $ system
by deriving the appropriate coupled f\/irst-order $q$-dif\/ference equations. We will utilise the form of the $E^{(1)}_6$
$q$-Painlev\'e system as given in \cite{KMNOY_2004} and \cite{KMNOY_2005} in terms of the variables~$f$,~$g$ under
the mapping
\begin{gather*}
( t, f,  g ) \mapsto \big( qt,  f(qt) \equiv \hat{f},  g(qt) \equiv \hat{g} \big) ,
\end{gather*}
and $ f(q^{-1}t) \equiv \check{f} $, etc. In these variables the coupled f\/irst-order
$q$-dif\/ference equations are
\begin{gather}
(g \check{f} -1) (gf -1)   = t^2\frac{(b_1g-1)(b_2g-1)(b_3g -1)(b_4g -1)}{(g-b_6t)(g-b_6^{-1}t)} ,
\label{eq:intro1} \\
(f \hat{g} -1) (fg -1)   = qt^2 \frac{(f- b_1)(f- b_2)(f -b_3)(f-b_4)}{(f- b_5qt)(f-b_5^{-1}t)} ,
\label{eq:intro2}
\end{gather}
with f\/ive independent parameters $b_1,\ldots,b_6$ subject to the constraint $b_1b_2b_3b_4=1$.

Our approach is to construct a sequence of $\tau$-functions starting with a deformation of a~specif\/ic weight
in the Askey table of hypergeometric orthogonal polynomial systems~\cite{KLS_2010}. However for the
purposes of the present work we will not explicitly exhibit these $\tau$-functions although one could do so
easily. The weight that we will take is the big $q$-Jacobi weight\footnote{However we will employ a dif\/ferent parameterisation of the
big $q$-Jacobi weight from that of the conventional form~(\ref{big-q-Jacobi}) in order that our results conform to the
the~$E^{(1)}_6$ $q$-Painlev\'e system as given by~(\ref{eq:intro1}),~(\ref{eq:intro2}); see~(\ref{DbqJ}).}
given by equation~(14.5.2) of~\cite{KLS_2010}
\begin{gather}
  w(x) = \frac{\big(a^{-1}x,c^{-1}x;q\big)_{\infty}}{\big(x,bc^{-1}x;q\big)_{\infty}} .
\label{big-q-Jacobi}
\end{gather}
The essential property of this weight, and the others in the Askey table, that we will utilise is that
they possess the $q$-analogue of the {\it semi-classical} property with respect to~$x$, namely that it
satisf\/ies the linear, f\/irst-order homogeneous $ q $-dif\/ference equation
\begin{gather*}
   \frac{w(qx)}{w(x)} = \frac{a(1-x)(c-bx)}{(a-x)(c-x)} ,
\end{gather*}
where the right-hand side is manifestly rational in $ x $.
Another feature of this weight is that the discrete lattice forming the support for the orthogonal polynomial
system is the $q$-linear lattice, one of four discrete quadratic lattices.
Consequently the perspective provided by our theoretical approach, then indicates that this case is the {\it master case} for
the $q$-linear lattices (as opposed to the $ D^{(1)}_5 $ system, for example) and all systems with such support will be degenerations of it.
The weight~(\ref{big-q-Jacobi}) has to be generalised, or {\it deformed}, in order to become relevant to $q$-Painlev\'e systems,
and such a generalisation turns out to introduce a new variable~$t$ and associated parameter so that it retains
the semi-classical character with respect to this variable. Using such a sequence of $\tau$-functions one employs
arguments to construct three systems of linear divided-dif\/ference equations which in turn characterise these.
One of these is the three-term recurrence relation of the polynomials orthogonal with respect to the deformed weight,
which in the Painlev\'e theory context is a distinguished Schlesinger transformation, while the two others are
our Lax pairs with respect to the spectral variable~$x$ and the deformation variable~$t$. Our method constructs a~specif\/ic
sequence of classical solutions to the~$ E^{(1)}_6 $ system and thus is technically valid for integer values of
a particular parameter, however we can simply analytically continue our results to the general case.

Lax pairs have been found for the $ E^{(1)}_6 $ system system using completely dif\/ferent techniques. In~\cite{Sa_2006} Sakai used a particular
degeneration of a two-variable Garnier extension to the Lax pairs for the $ D^{(1)}_5 $ $q$-Painlev\'e
system\footnote{This later system is also known as the $q$-$\rm P_{VI} $ system and its Lax pairs were constructed in~\cite{JS_1996}.}
(see~\cite{Sa_2005} for details on the multi-variable Garnier extension).
More recently Yamada \cite{Ya_2011} has reported Lax pairs for the~$ E^{(1)}_6 $ system by employing a degeneration starting
from a Lax pair for the
$ E^{(1)}_8 $ $q$-Painlev\'e equation through a sequence of limits~$ E^{(1)}_8 \to E^{(1)}_7 \to E^{(1)}_6 $.

The plan of our study is as follows. In Section~\ref{Theory} we recount the notations, def\/initions and basic facts of the general theory~\cite{Wi_2010a} in a self-contained manner omitting proofs. We draw heavily upon this theory in Section~\ref{Big_qJacobi} where we apply
it to the $q$-linear lattice and a natural extension or deformation of the big $q$-Jacobi weight. Again, using techniques
f\/irst expounded in~\cite{Wi_2010a}, we f\/ind explicit forms for the Lax pairs and verify the identif\/ication with the $E^{(1)}_6$
$q$-Painlev\'e system. At the conclusion of our study, in Section~\ref{Reconcile}, we relate our Lax pairs with those of both Sakai and Yamada.

\section{Deformed semi-classical OPS on quadratic lattices}\label{Theory}

We begin by summarising the key results of \cite{Wi_2010a}, in particular Sections~2,~3,~4 and~6 of that work,
which relate to semi-classical orthogonal polynomial systems with support on discrete, quadratic lattices.

Let $ \Pi_{n}[x] $ denote the linear space of polynomials in~$ x $ over~$ \C $ with
degree at most $ n\in\Z_{\geq 0} $. We def\/ine the {\it divided-difference operator} (DDO) $ \ddoAW_{x} $ by
\begin{gather}
   \ddoAW_{x} f(x) := \frac{f(\iota_{+}(x))-f(\iota_{-}(x))}{\iota_{+}(x)-\iota_{-}(x)} ,
\label{DDOdefn}
\end{gather}
and impose the condition that $ \ddoAW_{x}: \Pi_{n}[x] \to \Pi_{n-1}[x] $ for all $ n \in \N $.
In consequence we deduce that $ \iota_{\pm}(x) $ are the two $y$-roots of the quadratic equation
\begin{gather}
   \mathcal{A} y^2+2\mathcal{B} xy+\mathcal{C} x^2+2\mathcal{D} y+2\mathcal{E} x+\mathcal{F} = 0 .
\label{snul_Quad}
\end{gather}
Assuming $ \mathcal{A} \neq 0 $ the two $y$-roots $  y_{\pm} := \iota_{\pm}(x) $ for a given $ x $-value satisfy
\begin{gather*}
   \iota_{+}(x)+\iota_{-}(x) = -2\frac{\mathcal{B}x+\mathcal{D}}{\mathcal{A}} ,
 \qquad
   \iota_{+}(x)\iota_{-}(x) = \frac{\mathcal{C}x^2+2\mathcal{E}x+\mathcal{F}}{\mathcal{A}} ,
\end{gather*}
and their inverse functions $ \iota_{\pm}^{-1} $ are def\/ined by $ \iota_{\pm}^{-1}(\iota_{\pm}(x)) = x $.
For a given $ y $-value the quadra\-tic~(\ref{snul_Quad}) also def\/ines two $ x $-roots, if $ \mathcal{C}\neq 0 $,
which are consecutive points on the {\it $ x $-lattice},~$ x_s$,~$x_{s+1} $
parameterised by the variable $ s\in \Z $ and therefore def\/ines a map $ x_s \mapsto x_{s+1} $. Thus we have
the sequence of $x$-values $ \ldots ,x_{-2}, x_{-1}, x_0, x_1, x_2, \dots $ given by
$ \ldots $, $ \iota_{-}(x_0)=\iota_{+}(x_{-1}) $, $ \iota_{-}(x_1)=\iota_{+}(x_0) $, $ \ldots $
which we denote as the {\it lattice} or the {\it direct lattice} $ G $, and
the sequence of $y$-values $ \ldots ,y_{-2}, y_{-1}, y_0, y_1, y_2, \dots $ given by
$ \ldots $, $ y_0=\iota_{-}(x_0) $, $ y_1=\iota_{-}(x_1) $, $ y_2=\iota_{-}(x_2) $, $ \ldots $
as the {\it dual lattice} $ \tilde{G} $ (and distinct from the former in general).
A companion operator to the divided-dif\/ference operator $ \ddoAW_{x} $ is the mean or
average operator $ \moAW_{x} $ def\/ined by
\begin{gather*}
   \moAW_{x} f(x) = \tfrac{1}{2}\left[ f(\iota_{+}(x))+f(\iota_{-}(x)) \right] ,
\end{gather*}
so that the property $ \moAW_{x}: \Pi_{n}[x] \to \Pi_{n}[x] $ is ensured by the condition
we imposed upon $ \ddoAW_{x} $. The dif\/ference between consecutive points on the dual lattice
is given a distinguished notation through the def\/inition $ \Delta y(x) := \iota_{+}(x)-\iota_{-}(x) $.

We will also employ an operator notation for the mappings from points on the direct lattice to the
dual lattice $ E^{\pm}_{x} f(x):= f(\iota_{\pm}(x)) $ so that~\eqref{DDOdefn} can be written
\begin{gather*}
   \ddoAW_{x} f(x)
  = \frac{E^{+}_{x}f-E^{-}_{x}f}{E^{+}_{x}x-E^{-}_{x}x} ,
\end{gather*}
for arbitrary functions $ f(x) $.
The inverse functions $ \iota_{\pm}^{-1}(x) $ def\/ine operators $ (E^{\pm}_{x})^{-1} $ which map
points on the dual lattice to the direct lattice and also an adjoint to the
divided-dif\/ference operator $ \ddoAW_{x} $
\begin{gather*}
   \ddoAW_{x}^{*} f(x)
 := \frac{f\big(\iota_{+}^{-1}(x)\big)-f\big(\iota_{-}^{-1}(x)\big)}{\iota_{+}^{-1}(x)-\iota_{-}^{-1}(x)}
  = \frac{(E^{+}_{x})^{-1}f-(E^{-}_{x})^{-1}f}{(E^{+}_{x})^{-1}x-(E^{-}_{x})^{-1}x} .
\end{gather*}
The composite operators
$ E_{x}:=(E^{-}_{x})^{-1} E^{+}_{x} $ and $ E^{-1}_{x}=(E^{+}_{x})^{-1} E^{-}_{x} $
map between consecutive points on the direct lattice\footnote{However in the situation of a symmetric
quadratic $ \mathcal{A}=\mathcal{C} $ and $ \mathcal{D}=\mathcal{E} $, which
entails no loss of generality, then we have $ (E^{+}_{x})^{-1}=E^{-}_{x} $
and $ (E^{-}_{x})^{-1}=E^{+}_{x} $ and consequently there is no distinction
between the divided-dif\/ference operator and its adjoint.}.

Assuming $ \mathcal{A}\mathcal{C} \neq 0 $
one can classify these non-uniform quadratic lattices (or {\it SNUL}, special non-uniform lattices)
according to two parameters~-- the discriminant $ \mathcal{B}^2-\mathcal{A}\mathcal{C} $ and
\begin{gather*}
    \Theta = \det
               \begin{pmatrix} \mathcal{A} & \mathcal{B} & \mathcal{D} \\ \mathcal{B} & \mathcal{C} & \mathcal{E} \\ \mathcal{D} & \mathcal{E} & \mathcal{F}
               \end{pmatrix} ,
\end{gather*}
or $ \mathcal{A}\Theta = (\mathcal{B}^2-\mathcal{A}\mathcal{C})(\mathcal{D}^2-\mathcal{A}\mathcal{F})-(\mathcal{B}\mathcal{D}-\mathcal{A}\mathcal{E})^2 $.
The quadratic lattices are classif\/ied into four sub-cases \cite{Ma_1988,Ma_1995}:
$q$-quadratic ($ \mathcal{B}^2-\mathcal{A}\mathcal{C} \neq 0$ and $ \Theta < 0 $),
quadratic ($ \mathcal{B}^2-\mathcal{A}\mathcal{C} = 0$ and $ \Theta < 0 $),
$q$-linear ($ \mathcal{B}^2-\mathcal{A}\mathcal{C} \neq 0$ and $ \Theta = 0 $) and
linear ($ \mathcal{B}^2-\mathcal{A}\mathcal{C} = 0$ and $ \Theta = 0 $),
as the conic sections are divided into the elliptic/hyperbolic, parabolic, intersecting straight lines and parallel
straight lines respectively. The $q$-quadratic lattice, in its general non-symmetrical form, is the most general
case and the other lattices can be found from this by limiting processes.
For the quadratic class of lattices the parameterisation on $ s $ can be made explicit
using trigonometric/hyperbolic functions or their degenerations so we can employ
a parameterisation such that $ y_{s}=\iota_{-}(x_{s})=x_{s-1/2} $ and $ y_{s+1}=\iota_{+}(x_{s})=x_{s+1/2} $.
We denote the totality of lattice points
by $ G[x_0] := \{x_{s}:s\in \Z\} $ with the point $ x_{0} $ as the {\it basal point},
and of the dual lattice by $ \tilde{G}[x_0] := \{x_{s}:s\in \Z\+\frac{1}{2}\} $.

We def\/ine the {\it $ \ddoAW $-integral} of a
function def\/ined on the $x$-lattice $ f: G[x] \to \C $ with basal point~$ x_0 $ by the
Riemann sum over the lattice points
\begin{gather*}
   I[f](x_0) = \int_{G}  \ddoAW x\, f(x) := \sum_{s\in \Z} \Delta y(x_s)f(x_s) ,
\end{gather*}
where the sum is either a f\/inite subset of $ \Z $, namely $ \{ 0,\ldots,\Gt \} $, $ \Z_{\geq 0} $, or $ \Z $.
This def\/inition reduces to the usual def\/inition of the dif\/ference integral and
the Thomae--Jackson $q$-integrals in the canonical forms of the linear and
$q$-linear lattices respectively. Amongst a number of properties that f\/low from this def\/inition we have
an analog of the fundamental theorem of calculus
\begin{gather}
   \int_{x_0\leq x_{s} \leq x_{\Gt}} \ddoAW x\, \ddoAW_{x} f(x) = f(E^{+}_{x}x_{\Gt})-f(E^{-}_{x}x_{0}) .
\end{gather}

Central to our study are orthogonal polynomial systems (OPS) def\/ined on $G$, and a general reference
for a background on these and other considerations is the monograph by Ismail \cite{Ismail_2005}.
Our OPS is def\/ined via orthogonality relations with support on $ G $
\begin{gather*}
  \int_{G} \ddoAW x\; w(x)p_n(x)\,l_m(x) = \begin{cases} 0, \  & 0\leq m<n, \\ h_n, \  & m=n, \end{cases}
  \qquad n \geq 0, \qquad h_n \neq 0,
\end{gather*}
where $ \{l_m(x)\}^{\infty}_{m=0} $ is any system of polynomial bases with exact $ \operatorname{deg}_{x}(l_m)=m $.
Such relations def\/ine a sequence of orthogonal polynomials $ \{p_n(x)\}^{\infty}_{n=0} $ under suitable conditions (see~\cite{Ismail_2005}).
An immediate consequence of orthogonality is that the orthogonal polynomials satisfy a three term recurrence relation of the form
\begin{gather}
   a_{n+1}p_{n+1}(x) = (x-b_n)p_n(x) - a_np_{n-1}(x), \qquad n \geq 0,\nonumber\\
    a_n \neq 0 , \qquad p_{-1} = 0, \qquad p_0 = \gamma_0 .
\label{3termRR}
\end{gather}
However we require non-polynomial solutions to this linear second-order dif\/ference equation, which are
linearly independent of the polynomial solutions. To this end we def\/ine the {\it Stieltjes function}
\begin{gather*}
  f(x) \equiv \int_{G} \ddoAW y\, \frac{w(y)}{x-y} , \qquad x \notin G .
\end{gather*}
A set of non-polynomial solutions to (\ref{3termRR}), termed {\it associated functions} or {\it functions of the second kind},
and which generalise the Stieltjes function, are given by
\begin{gather*}
  q_n(x) \equiv \int_{G} \ddoAW y\, w(y)\frac{p_n(y)}{x-y} , \qquad n \geq 0, \qquad x \notin G .
\end{gather*}
The associated function solutions dif\/fer from the orthogonal polynomial solutions in that they have
the initial conditions $ q_{-1} = 1/a_{0}\gamma_{0} $, $ q_{0} = \gamma_{0}f $.
The utility and importance of the Stieltjes function lies in the fact that that it connects $ p_{n} $ and
$ q_{n} $ whereby the dif\/ference $ fp_{n} - q_{n} $ is exactly a polynomial of degree $ n-1 $ which
itself satisf\/ies ~(\ref{3termRR}) in place of $ p_{n} $. This relation is crucial for the arguments adopted in~\cite{Wi_2010a}.
With the polynomial and non-polynomial solutions we form the $ 2\times 2 $ matrix variable, which occupies a primary
position in our theory:
\begin{gather*}
   Y_n(x) = \begin{pmatrix} p_n(x) & \dfrac{q_{n}(x)}{w(x)} \\
                              p_{n-1}(x) & \dfrac{q_{n-1}(x)}{w(x)}
            \end{pmatrix} .
\end{gather*}
In this matrix variable the three-term recurrence relation takes the form
\begin{gather}
    Y_{n+1} = K_{n}Y_{n}, \qquad
     K_n(x) = \frac{1}{a_{n+1}}
              \begin{pmatrix} x-b_n & -a_{n} \\
                              a_{n+1} & 0
              \end{pmatrix}, \qquad \det K_n = \frac{a_n}{a_{n+1}} .
\label{RR}
\end{gather}
A key result required in the analysis of OPS are the expansions of polynomial solutions about the f\/ixed singularity at $ x=\infty $
\begin{gather}
   p_n(x) = \gamma_n \left[ x^{n} - \left( \sum^{n-1}_{i=0}b_i \right)x^{n-1}
   + \left( \sum_{0\leq i<j<n}b_ib_j-\sum^{n-1}_{i=1}a^2_i \right)x^{n-2}
   + {\rm O}\big(x^{n-3}\big) \right] ,
\label{ops_pExp}
\end{gather}
valid for $ n \geq 1 $, while for the associated functions the expansions read
\begin{gather}
   q_n(x) = \gamma^{-1}_n \left[ x^{-n-1} + \!\left( \sum^{n}_{i=0}b_i \!\right)x^{-n-2}
   + \!\left(\! \sum_{0\leq i\leq j\leq n}\!\!\!b_ib_j+\sum^{n+1}_{i=1}a^2_i \!\right)x^{-n-3}
   + {\rm O}\big(x^{-n-4}\big) \!\right] ,\!\!\!
\label{ops_eExp}
\end{gather}
valid for $ n \geq 0 $.

In order to proceed any further we need to impose some structure on the weight characterising our OPS~--
in particular its spectral characteristics~-- and this takes the form of the def\/inition of a~{\it $\ddoAW$-semi-classical weight}~\cite{Ma_1988}. Such a weight satisf\/ies a f\/irst-order homogeneous
divided-dif\/ference equation
\begin{gather}
   \frac{w(y_{+})}{w(y_{-})} = \frac{W+\Delta yV}{W-\Delta yV}(x) ,
\label{sc_Defn}
\end{gather}
where $ W(x)$, $V(x)$ are irreducible polynomials in the spectral variable $ x $, which we call
{\it spectral polynomials}. As a consequence of this,
under reasonable assumptions on the parameters of the weight, the Stieltjes function satisf\/ies
an inhomogeneous form of (\ref{sc_Defn})
\begin{gather}
   W\ddoAW_{x} f = 2V\moAW_{x} f + U ,
\label{sc_Stieltjes}
\end{gather}
where in addition $ U(x) $ is a polynomial of~$ x $.
A {\it generic} or {\it regular $\ddoAW$-semi-classical weight} has two properties:
\begin{itemize}\itemsep=0pt
\item[(i)]
  strict inequalities in the degrees of the spectral data polynomials, i.e.,
  $ \deg_x W = M $, $ \deg_x V = M-1 $ and $\deg_x U = M-2 $, and
\item[(ii)]
  the lattice generated by any zero of $ (W^2-\Delta y^2V^2)(x) $, say $ \tilde{x}_1 $, does not coincide
  with another zero, $ \tilde{x}_2 $, i.e.\ if $ (W^2-\Delta y^2V^2)(\tilde{x}_2)=0 $ then
  $ \tilde{x}_2 \notin \iota^{2\mathbb{Z}}_{\pm}\tilde{x}_1 $.
\end{itemize}

Further consequences of semi-classical assumptions are a system of spectral divided-dif\/ference
equations for the matrix variable $ Y_{n} $, i.e., the {\it spectral divided-difference} equation
\begin{gather}
   \ddoAW_{x} Y_n(x) := A_{n}\moAW_{x} Y_n(x)\nonumber\\
\hphantom{\ddoAW_{x} Y_n(x)}{}
   = \frac{1}{W_{n}(x)}
     \begin{pmatrix} \Omega_n(x) & -a_n\Theta_n(x) \\
                     a_n\Theta_{n-1}(x) & -\Omega_n(x)-2V(x)
     \end{pmatrix}\moAW_{x} Y_n(x),
  \qquad n \geq 0 ,
\label{spectral_DDO}
\end{gather}
with $ A_{n} $ termed the {\it spectral matrix}.
For the $\ddoAW$-semi-classical class of weights the coef\/f\/icients appearing in the spectral matrix,
$ W_{n}$, $\Omega_n$, $\Theta_n $, are polynomials in~$ x $, with f\/ixed degrees independent of the index~$ n $.
These {\it spectral coefficients} have terminating expansions about $ x=\infty $ with the
leading order terms
\begin{gather}
  W_n(x) = \tfrac{1}{2}W + \tfrac{1}{4}[W+\Delta yV]\left( \frac{y_{+}}{y_{-}} \right)^{n}\nonumber\\
\hphantom{W_n(x) =}{}
                     + \tfrac{1}{4}[W-\Delta yV]\left( \frac{y_{-}}{y_{+}} \right)^{n}
                     + {\rm O}\big(x^{M-1}\big),
  \qquad n \geq 0 ,
\label{spectral_coeffEXP:a}
\\
  \Theta_n(x) =   \frac{1}{y_{-}\Delta y}[W+\Delta yV]\left( \frac{y_{+}}{y_{-}} \right)^{n}\nonumber\\
\hphantom{\Theta_n(x) =}{}
             - \frac{1}{y_{+}\Delta y}[W-\Delta yV]\left( \frac{y_{-}}{y_{+}} \right)^{n}
             + {\rm O}\big(x^{M-3}\big),
  \qquad n \geq 0 ,
\label{spectral_coeffEXP:b}
\\
  \Omega_n(x)+V(x) =   \frac{1}{2\Delta y}[W+\Delta yV]\left( \frac{y_{+}}{y_{-}} \right)^{n}\nonumber\\
\hphantom{\Omega_n(x)+V(x) =}{}
               - \frac{1}{2\Delta y}[W-\Delta yV]\left( \frac{y_{-}}{y_{+}} \right)^{n}
               + {\rm O}\big(x^{M-2}\big),
  \qquad n \geq 0 ,
\label{spectral_coeffEXP:c}
\end{gather}
where $ M = \deg_x(W_{n}) $.

Compatibility of the spectral divided-dif\/ference equations (\ref{spectral_DDO}) and recurrence rela\-tions~(\ref{RR}) imply that the spectral matrix and the recurrence matrix satisfy
\begin{gather}
  K_n(y_+)\left( 1-\tfrac{1}{2}\Delta y A_{n} \right)^{-1}\left( 1+\tfrac{1}{2}\Delta y A_{n} \right)\nonumber\\
  \qquad{}=
  \left( 1-\tfrac{1}{2}\Delta y A_{n+1} \right)^{-1}\left( 1+\tfrac{1}{2}\Delta y A_{n+1} \right)K_n(y_-) ,
  \qquad n \geq 0 .
\label{spectralA_recur:b}
\end{gather}
These relations can be rewritten in terms of the spectral coef\/f\/icients arising in \eqref{spectral_DDO} as
recurrence relations in $ n $,
\begin{gather}
  W_{n+1} = W_{n}+\tfrac{1}{4}\Delta y^2 \Theta_{n},
  \qquad n \geq 0 ,
\label{spectral_coeff_recur:a}  \\
  \Omega_{n+1}+\Omega_{n}+2V = (\moAW_{x} x-b_{n})\Theta_{n},
  \qquad n \geq 0 ,
\label{spectral_coeff_recur:b}  \\
  (W_{n}\Omega_{n+1}-W_{n+1}\Omega_{n})(\moAW_{x} x-b_{n})
\nonumber\\ \qquad
 {} = - \tfrac{1}{4}\Delta y^2 \Omega_{n+1}\Omega_{n}
  + W_{n}W_{n+1} + a^2_{n+1}W_{n}\Theta_{n+1} - a^2_{n}W_{n+1}\Theta_{n-1},
  \qquad n \geq 0 .
\label{spectral_coeff_recur:c}
\end{gather}
Another important deduction from these relations is that the spectral coef\/f\/icients satisfy a~bilinear
relation
\begin{gather}
 W_{n}(W_{n}-W)
  = -\tfrac{1}{4}\Delta y^2\det
     \begin{pmatrix} \Omega_n & -a_n\Theta_n \\
                     a_n\Theta_{n-1} & -\Omega_n-2V
     \end{pmatrix},
  \qquad n \geq 0 .
\label{spectral_Det}
\end{gather}
The matrix product appearing in (\ref{spectralA_recur:b}), and recurring subsequently, is called the
{\it Cayley transform} of~$ A_{n} $ and it has the evaluation
\begin{gather}
  \left( 1-\tfrac{1}{2}\Delta y A_{n} \right)^{-1}\left( 1+\tfrac{1}{2}\Delta y A_{n} \right)
 \label{spectral_prod}\\
\qquad{}  = \frac{1}{W+\Delta y V}
     \begin{pmatrix} 2W_n-W+\Delta y (\Omega_n+V) & -\Delta y  a_n\Theta_n \\
                     \Delta y a_n\Theta_{n-1} & 2W_n-W-\Delta y (\Omega_n+V)
     \end{pmatrix},
  \qquad n \geq 0 .\nonumber
\end{gather}
This result motivates the following def\/initions
\begin{gather}
  \mathfrak{W}_{\pm} := 2W_n-W\pm\Delta y (\Omega_n+V), \qquad
  \mathfrak{T}_{+} := \Delta y  a_n\Theta_n, \nonumber\\
  \mathfrak{T}_{-} := \Delta y  a_n\Theta_{n-1},
  \qquad n \geq 1 ,
\label{spectral_def}
\end{gather}
whilst for $ n=0 $ we have
$ \mathfrak{W}_{\pm}(n=0) := W\pm\Delta y V $,
$ \mathfrak{T}_{+}(n=0) := -\Delta y a_0\gamma^2_0 U $, and
\mbox{$ \mathfrak{T}_{-}(n=0) := 0 $}. Thus we def\/ine
\begin{gather}
  A^{*}_n :=
     \begin{pmatrix} \mathfrak{W}_{+} & -\mathfrak{T}_{+} \\
                     \mathfrak{T}_{-} &  \mathfrak{W}_{-}
     \end{pmatrix} .
\label{AC_defn}
\end{gather}
In a scalar formulation of the matrix linear divided-dif\/ference equation \eqref{spectral_DDO} one of the
components, $p_n$ say, satisf\/ies a linear second-order divided-dif\/ference equation
of the form
\begin{gather}
  E^{+}_{x} \left( \frac{W + \Delta y V}{\Delta y \Theta_n} \right) (E^{+}_{x})^2 p_n +
	E^{-}_{x} \left( \frac{W - \Delta y V}{\Delta y \Theta_n} \right) (E^{-}_{x})^2 p_n\nonumber\\
\qquad{}
    - \left\{ E^{+}_{x} \left( \frac{ \mathfrak{W}_+}{\Delta y \Theta_n} \right)
	+ E^{-}_{x} \left( \frac{ \mathfrak{W}_-}{\Delta y \Theta_n} \right)
      \right\} E^{+}_{x}E^{-}_{x} p_n = 0.
\label{eq:LSODDE}
\end{gather}

Thus far our theoretical construction can only account for the OPS occurring in the Askey table~--
the hypergeometric and basic hypergeometric orthogonal polynomial systems~\cite{KLS_2010}. To
step beyond these, and in particular to make contact with the discrete Painlev\'e systems, one
has to introduce pairs of deformation variables and parameters into the OPS.
We denote such a single deformation variable by~$ t $, def\/ined on a quadratic lattice (and possibly distinct from that of
the spectral variable), with advanced and retarded nodes at $ \iota_{\pm}(t) = u_{\pm} $, $ \Delta u = \iota_{+}(t)-\iota_{-}(t) $.
We introduce such deformations with imposed structures that are analogous to those of the spectral variable.
Thus, corresponding to the def\/inition~(\ref{sc_Defn}), we deem that a
{\it deformed  $\ddoAW$-semi-classical weight}~$w(x;t) $ satisf\/ies the additional f\/irst-order homogeneous
divided-dif\/ference equation
\begin{gather}
  \frac{w(x;u_{+})}{w(x;u_{-})} = \frac{R+\Delta uS}{R-\Delta uS}(x;t) ,
\label{dsc_Defn}
\end{gather}
where the deformation data polynomials, $ R(x;t)$, $S(x;t) $, are irreducible polynomials in~$ x $.
The spectral data polynomials, $ W(x;t)$, $V(x;t) $, and the deformation data polynomials,
$ R(x;t)$, $S(x;t) $, now must satisfy the compatibility relation
\begin{gather}
  \frac{W+\Delta yV}{W-\Delta yV}(x;u_{+})
  \frac{R+\Delta uS}{R-\Delta uS}(y_{-};t)
  =
  \frac{W+\Delta yV}{W-\Delta yV}(x;u_{-})
  \frac{R+\Delta uS}{R-\Delta uS}(y_{+};t) .
\label{WV-RS}
\end{gather}
The deformed $\ddoAW$-semi-classical deformation condition that corresponds to (\ref{sc_Stieltjes}) is that the
Stieltjes transform satisf\/ies an inhomogeneous version of (\ref{dsc_Defn})
\begin{gather*}
    R\ddoAW_{t} f = 2S\moAW_{t} f + T ,
%\label{dsc_Stieltjes}
\end{gather*}
with $ T(x;t) $ being an irreducible polynomial in $ x $ with respect to $ R $ and $ S $.
Compatibility of spectral and deformation divided-dif\/ference equations for $ f $ implies the following
identity on $ U $ and $ T $
\begin{gather*}
     \Delta y\left[ \frac{(W+\Delta yV)(x;u_{+})}{(W+\Delta yV)(x;u_{-})}(R+\Delta uS)(y_{-};t)U(x;u_{-})
                  -(R-\Delta uS)(y_{-};t)U(x;u_{+}) \right]
 \\
\qquad{} = \Delta u\left[ (W+\Delta yV)(x;u_{+})T(y_{-};t)
                  -(W-\Delta yV)(x;u_{+})\frac{(R-\Delta uS)(y_{-};t)}{(R-\Delta uS)(y_{+};t)}T(y_{+};t) \right] .
\end{gather*}
Corresponding to the (\ref{spectral_DDO}) the
deformed $\ddoAW$-semi-classical OPS satisf\/ies the {\it deformation divided-difference} equation
\begin{gather}
   \ddoAW_{t} Y_n := B_{n}\moAW_{t} Y_n
   = \frac{1}{R_{n}}
     \begin{pmatrix} \Gamma_n & \Phi_n \\
                     \Psi_{n} & \Xi_n
     \end{pmatrix}\moAW_{t} Y_n ,
  \qquad n \geq 0 .
\label{deform_DDO}
\end{gather}
The deformation coef\/f\/icients appearing in matrix $B_n$ above satisfy a linear identity
\begin{gather}
   \Psi_{n} = -\frac{a_{n}}{a_{n-1}}\Phi_{n-1},
  \qquad n \geq 1 ,
\label{Defm_offD}
\end{gather}
and a trace identity
\begin{gather*}
  \Delta u (\Gamma_n+\Xi_n)
 = 2H_n\left[
       \frac{R+\Delta uS}{a_n(u_{-})}-\frac{R-\Delta uS}{a_n(u_{+})} \right],
  \qquad n \geq 0 ,
%\label{Defm_Tr}
\end{gather*}
which means that only three of these are independent. Here $ H_{n} $ is a constant with respect to $ x $
and arises as a decoupling constant which will be set subsequently in applications to a convenient value.
The deformation coef\/f\/icients are all
polynomials in $ x $, with f\/ixed degrees independent of the index $ n $ but with non-trivial
$ t $ dependence. Let $ L = \max(\deg_x R,\deg_x S) $.
As $ x \to \infty $ we have the leading orders of the terminating expansions of
the following deformation coef\/f\/icients
\begin{gather}
  \frac{2}{H_n}R_n =
  -(\gamma_n(u_{+})+\gamma_n(u_{-}))\left[ \frac{R-\Delta uS}{\gamma_{n-1}(u_{+})}+\frac{R+\Delta uS}{\gamma_{n-1}(u_{-})}
                                \right]
  + {\rm O}\big(x^{L-1}\big)
  \qquad n \geq 0 ,
\label{deform_coeff_Exp:a}
\\
  \frac{\Delta u}{2H_n}\Phi_n =
   \left[ (R+\Delta uS)\frac{\gamma_{n}(u_{+})}{\gamma_{n}(u_{-})}-(R-\Delta uS)\frac{\gamma_{n}(u_{-})}{\gamma_{n}(u_{+})}
   \right]x^{-1}
  + {\rm O}\big(x^{L-2}\big),
  \qquad n \geq 0 ,
\label{deform_coeff_Exp:b}
\end{gather}
and
\begin{gather}
  \frac{\Delta u}{H_n}\Gamma_n =
   (\gamma_n(u_{-})-\gamma_n(u_{+}))\left[ \frac{R+\Delta uS}{\gamma_{n-1}(u_{-})}+\frac{R-\Delta uS}{\gamma_{n-1}(u_{+})}
                                \right]
  + {\rm O}\big(x^{L-1}\big),
  \qquad n \geq 0 .
\label{deform_coeff_Exp:d}
\end{gather}
Compatibility of the deformation divided-dif\/ference equation (\ref{deform_DDO}) and the recurrence relation~(\ref{RR}) implies the relation
\begin{gather}
  K_n(;u_{+})\left( 1-\tfrac{1}{2}\Delta uB_{n} \right)^{-1}\left( 1+\tfrac{1}{2}\Delta uB_{n} \right)\nonumber\\
  \qquad{}
  =
  \left( 1-\tfrac{1}{2}\Delta uB_{n+1} \right)^{-1}\left( 1+\tfrac{1}{2}\Delta uB_{n+1} \right)K_n(;u_{-}),
  \qquad n \geq 0 .
\label{deform_Brecur:b}
\end{gather}
From this we can deduce that the deformation coef\/f\/icients, $ R_n$, $\Gamma_n$, $\Phi_n $, satisfy recurrence
relations in $ n $ in parallel to those of~(\ref{spectral_coeff_recur:a}),~(\ref{spectral_coeff_recur:b})
\begin{gather*}
  \frac{a_{n+1}(u_{-})}{H_{n+1}}(-2R_{n+1}+\Delta u\Gamma_{n+1})
 +\frac{a_{n}(u_{-})}{H_{n}}(2R_{n}+\Delta u\Gamma_{n})
\nonumber \\
\qquad{} = -[x-b_n(u_{-})]\frac{\Delta u}{H_n}\Phi_n
   +2a_n(u_{-})\left( \frac{R+\Delta uS}{a_n(u_{-})}-\frac{R-\Delta uS}{a_n(u_{+})} \right) ,
  \qquad n \geq 0 ,
%\label{Defm_recur:a}
\\
  \frac{a_{n+1}(u_{+})}{H_{n+1}}(2R_{n+1}+\Delta u\Gamma_{n+1})
 +\frac{a_{n}(u_{+})}{H_{n}}(-2R_{n}+\Delta u\Gamma_{n})
 \nonumber\\
\qquad{} = -[x-b_n(u_{+})]\frac{\Delta u}{H_n}\Phi_n
   +2a_n(u_{+})\left( \frac{R+\Delta uS}{a_n(u_{-})}-\frac{R-\Delta uS}{a_n(u_{+})} \right) ,
  \qquad n \geq 0 .
%\label{Defm_recur:b}
\end{gather*}
The deformation coef\/f\/icients satisfy the bilinear or determinantal identity
\begin{gather*}
 R^2_n+\tfrac{1}{4}\Delta u^2 \left[ \Gamma_{n}\Xi_{n}-\Phi_{n}\Psi_{n} \right]
 = - H_nR_n\left[
       \frac{R+\Delta uS}{a_n(u_{-})}+\frac{R-\Delta uS}{a_n(u_{+})} \right],
  \qquad n \geq 0 ,
%\label{deform_Det}
\end{gather*}
which is the analogue of (\ref{spectral_Det}). The matrix product given in (\ref{deform_Brecur:b})
has the evaluation
\begin{gather*}
  \left( 1-\tfrac{1}{2}\Delta u B_{n} \right)^{-1}\left( 1+\tfrac{1}{2}\Delta u B_{n} \right)
  = \frac{a_n(u_{-})}{2H_n(R+\Delta u S)}  \\
 \qquad\times    \begin{pmatrix} 2R_n+2H_n\dfrac{R-\Delta u S}{a_n(u_{+})}+\Delta u \Gamma_n & \Delta u \Phi_n \\
                     \Delta u \Psi_n & 2R_n+2H_n\dfrac{R+\Delta u S}{a_n(u_{-})}-\Delta u \Gamma_n .
     \end{pmatrix},
  \qquad n \geq 0 .
%\label{deform_prod}
\end{gather*}
This again motivates the def\/initions
\begin{gather}
  \mathfrak{R}_{\pm} := 2R_n+2H_n\frac{R\mp\Delta u S}{a_n(u_{\pm})}\pm\Delta u \Gamma_n, \nonumber\\
  \mathfrak{P}_{+} :=-\Delta u \Phi_n, \qquad
  \mathfrak{P}_{-} := \Delta u \Psi_n, \qquad n \geq 1,
\label{deform_def:abc}
\end{gather}
together with
\begin{gather*}
  B^{*}_n :=
     \begin{pmatrix} \mathfrak{R}_{+} & -\mathfrak{P}_{+} \\
                     \mathfrak{P}_{-} &  \mathfrak{R}_{-}
     \end{pmatrix} .
%\label{BC_defn}
\end{gather*}

Our f\/inal relation expresses the compatibility of the spectral and deformation divided-dif\/fe\-ren\-ce equations.
The spectral matrix $ A_n(x;t) $ and the deformation matrix $ B_n(x;t) $ satisfy the
$\ddoAW$-{\it Schlesinger equation}
\begin{gather}
  \left( 1-\tfrac{1}{2}\Delta yA_{n}(;u_{+}) \right)^{-1}\left( 1+\tfrac{1}{2}\Delta yA_{n}(;u_{+}) \right)
  \left( 1-\tfrac{1}{2}\Delta uB_{n}(y_{-};) \right)^{-1}\left( 1+\tfrac{1}{2}\Delta uB_{n}(y_{-};) \right) \label{spectral+deform:b}\\
  =
  \left( 1-\tfrac{1}{2}\Delta uB_{n}(y_{+};) \right)^{-1}\left( 1+\tfrac{1}{2}\Delta uB_{n}(y_{+};) \right)
  \left( 1-\tfrac{1}{2}\Delta yA_{n}(;u_{-}) \right)^{-1}\left( 1+\tfrac{1}{2}\Delta yA_{n}(;u_{-}) \right) .\nonumber
\end{gather}
Let us def\/ine the quotient
\begin{gather*}
  \chi
    \equiv \frac{(W+\Delta yV)(x;u_{+})}{(W+\Delta yV)(x;u_{-})}
           \frac{(R+\Delta u S)(y_{-};t)}{(R+\Delta u S)(y_{+};t)}
    = \frac{(W-\Delta yV)(x;u_{+})}{(W-\Delta yV)(x;u_{-})}
      \frac{(R-\Delta u S)(y_{-};t)}{(R-\Delta u S)(y_{+};t)} .
%\label{twist}
\end{gather*}
The compatibility relation (\ref{spectral+deform:b}) can be rewritten as the matrix equation
\begin{gather}
 \chi  B^*_n(y_{+};t)A^*_n(x;u_{-}) = A^*_n(x;u_{+})B^*_n(y_{-};t) ,
\label{A-B}
\end{gather}
or component-wise with the new variables in the more practical form as
\begin{gather}
 \chi \left[
 \mathfrak{W}_{+}(x;u_{-})\mathfrak{R}_{+}(y_{+};t) - \mathfrak{T}_{-}(x;u_{-})\mathfrak{P}_{+}(y_{+};t)
           \right]\nonumber\\
\qquad{}
 = \mathfrak{W}_{+}(x;u_{+})\mathfrak{R}_{+}(y_{-};t) - \mathfrak{T}_{+}(x;u_{+})\mathfrak{P}_{-}(y_{-};t) ,
\label{S+D:a}
\\
 \chi \left[
 \mathfrak{T}_{+}(x;u_{-})\mathfrak{R}_{+}(y_{+};t) + \mathfrak{W}_{-}(x;u_{-})\mathfrak{P}_{+}(y_{+};t)
           \right]\nonumber\\
\qquad{}
 = \mathfrak{T}_{+}(x;u_{+})\mathfrak{R}_{-}(y_{-};t) + \mathfrak{W}_{+}(x;u_{+})\mathfrak{P}_{+}(y_{-};t) ,
\label{S+D:b}
\\
 \chi \left[
 \mathfrak{T}_{-}(x;u_{-})\mathfrak{R}_{-}(y_{+};t) + \mathfrak{W}_{+}(x;u_{-})\mathfrak{P}_{-}(y_{+};t)
           \right]\nonumber\\
\qquad{}
 = \mathfrak{T}_{-}(x;u_{+})\mathfrak{R}_{+}(y_{-};t) + \mathfrak{W}_{-}(x;u_{+})\mathfrak{P}_{-}(y_{-};t) ,
\label{S+D:c}
\\
 \chi \left[
 \mathfrak{W}_{-}(x;u_{-})\mathfrak{R}_{-}(y_{+};t) - \mathfrak{T}_{+}(x;u_{-})\mathfrak{P}_{-}(y_{+};t)
           \right]\nonumber\\
\qquad{}
 = \mathfrak{W}_{-}(x;u_{+})\mathfrak{R}_{-}(y_{-};t) - \mathfrak{T}_{-}(x;u_{+})\mathfrak{P}_{+}(y_{-};t) .
\label{S+D:d}
\end{gather}

For a general quadratic lattice there exists two f\/ixed points def\/ined by $ \iota_{+}(x)=\iota_{-}(x) $,
and let us denote these two points of the $x$-lattice by $ x_L $ and $ x_R $.
By analogy with the linear lattices we conjecture the existence of
fundamental solutions to the spectral divided-dif\/ference equation about $ x=x_{L},x_{R} $ which we
denote by $ Y_{L}$, $Y_{R} $ respectively.
Furthermore let us def\/ine the {\it connection matrix}
\begin{gather*}
  P(x;t) := Y_{R}(x;t)^{-1}Y_{L}(x;t) .
\end{gather*}
From the spectral divided-dif\/ference equation (\ref{spectral_DDO}) it is clear that $ P $ is a $\ddoAW$-constant function with respect
to $ x $, that is to say
\begin{gather*}
  P(y_{+};t) = P(y_{-};t) .
\end{gather*}
In addition it is clear from the deformation divided-dif\/ference equation~(\ref{deform_DDO}) that this type of deformation is also a
{\it connection preserving deformation} in the sense that
\begin{gather*}
  P(x;u_{+}) = P(x;u_{-}) .
\end{gather*}
This is our analogue of the monodromy matrix and generalises the connection matrix of Birkhof\/f and his
school~\cite{Birkhoff_1911, Birkhoff_1913}, although we emphasise that we have made an empirical
observation of this fact and not provided any rigorous statement of it.

\section[Big $q$-Jacobi OPS]{Big $\boldsymbol{q}$-Jacobi OPS}\label{Big_qJacobi}

As our central reference on the Askey table of basic hypergeometric orthogonal polynomial systems we employ~\cite{KS_1998}, or its
modern version~\cite{KLS_2010}. We consider a sub-case of the quadratic lattices, in particular the
$q$-linear lattice in both the spectral and deformation variables~$ x $ and~$ t $ in its standardised form, so that
$ \iota_{+}(x)=qx$, $\iota_{-}(x)=x$, $\Delta y(x)=(q-1)x $ and $ \iota_{+}(t)=qt$, $\iota_{-}(t)=t$, $\Delta u(t)=(q-1)t $.
In \cite{KLS_2010} the big $q$-Jacobi weight given by equation~(14.5.2) is
\begin{gather*}
  w(x) = \frac{(a^{-1}x,c^{-1}x;q)_{\infty}}{(x,bc^{-1}x;q)_{\infty}} ,
\end{gather*}
subject to $ 0 < aq,bq < 1 $, $ c < 0 $ with respect to the Thomae--Jackson $q$-integral
\begin{gather*}
 \int^{aq}_{bq} d_{q}x \, f(x) .
\end{gather*}
The $q$-shifted factorials have the standard def\/inition
\begin{gather*}
  (a;q)_{\infty} = \prod^{\infty}_{j=0}\big(1-aq^j\big), \qquad |q|<1, \qquad (a_1,\ldots,a_n;q)_{\infty} = (a_1;q)_{\infty}\cdots (a_n;q)_{\infty} .
\end{gather*}
We deform this weight by introducing an extra $q$-shifted factorial into the numerator and denominator containing
the deformation variable and parameter, and relabeling the big $q$-Jacobi parameters. We propose the following weight
\begin{gather}
 w(x;t)=\frac{\big(b_2 x,b_3 x,b_6^{-1}xt^{-1};q\big)_{\infty}}{\big(b_1 x,b_4 x,b_6 xt^{-1};q\big)_{\infty}} .
\label{DbqJ}
\end{gather}
A condition $ b_1b_2b_3b_4 = 1 $ will apply, so we have four free parameters.
We do not need to specify the support for this weight for the purposes of our work, but suf\/f\/ice it to
say that any Thomae--Jackson $q$-integral with terminals coinciding with any pair of zeros and poles
of the weight would be suitable.

The spectral data polynomials are computed to be
\begin{gather}
 W+\Delta yV  = b_6\left(1-b_1x\right)\left(1-b_4 x\right)\left(t- b_6x\right) ,
%\label{eq:weight2}
\nonumber \\
 W-\Delta yV  = \left(1-b_2x\right)\left(1-b_3x\right)\left(b_6t-x\right) .
\label{eq:weight3}
\end{gather}
Clearly the regular $ M=3 $ case is applicable and we seek solutions to the spectral coef\/f\/icients with
$\deg_x W_n = 3$, $\deg_x \Omega_n = 2$, $\deg_x \Theta_n = 1 $. Our procedure is to employ
the following algorithm, as detailed in~\cite{Wi_2010a}. Firstly we parameterise the
spectral matrix in a minimal way; secondly we relate the parameterisation of the deformation matrix to
that of the spectral matrix and thus close the system of unknowns; and f\/inally utilise these
parameterisations in the system of over-determined equations to derive evolution equations for our
primary variables. What constitutes the primary variables will emerge from the calculations themselves.

\begin{proposition}
Let us define a new parameter $ b_5 $ replacing $q^n$ by
\begin{gather*}
 q^n = \frac{b_5}{b_1b_4b_6}, \qquad n \in \Z_{\geq 0} .
\end{gather*}
Let the parameters satisfy the conditions $ q \neq 1 $, $ b_5 \neq q^{-1/2},\pm 1,q^{1/2} $,
$ b_1b_4 \neq 0, \infty $ and $ b_2b_3 \neq 0, \infty $.
Given the degrees of the spectral coefficients we parameterise these by
\begin{gather*}
   2W_{n}-W = w_3x^3+w_2x^2+w_1x+w_0 ,
  \\
   \Omega_{n}+V = v_2x^2+v_1x+v_0 ,
  \\
   \Theta_{n} = u_{1}(x-\lambda_n) .
\end{gather*}
Let $ \lambda_n $ be the unique zero of the $ (1,2) $ component of $ A^*_n $, i.e., $ \Theta_n(x) $ and define
the further variables $ \nu_n = (2W_{n}-W)(\lambda_{n},t) $ and $ \mu_{n} =(\Omega_{n}+V)(\lambda_{n},t) $.
Then the spectral coefficients are given by
\begin{gather}
  2W_{n}-W = x^2 \frac{\nu_n}{\lambda_n^2}
  +\tfrac{1}{2}(x-\lambda_n)\nonumber\\
 \hphantom{2W_{n}-W =}{}
  \times\left[ -b_6\big(b_5+b_5^{-1}\big)x^2+\frac{1+\left(b_1+b_2+b_3+b_4\right)b_6t+b_6^2}{\lambda_n}x
    -2t b_6\frac{x+\lambda_n}{\lambda_n^2} \right] ,
\label{Wparam}
\\
 \Omega_{n}+V = \mu_n
  +\frac{b_6}{2b_5 (1-q)\big(1-b_5^2\big)\lambda_n^2}\left(x-\lambda_n\right) \Big\{ {-}\big(1-b_5^2\big)^2\lambda_n^2x
\nonumber\\ \phantom{\Omega_{n}+V =}{}
    -2 b_5^2\left[ b_1^{-1}+b_2^{-1}+b_3^{-1}+b_4^{-1}+\big(b_6+b_6^{-1}\big)t-2\lambda_n \right]\lambda_n^2
\nonumber\\
\phantom{\Omega_{n}+V =}{}
      +b_5b_6^{-1}\big(1+b_5^2\big)\left[  \big(1+ (b_1+b_2+b_3+b_4 )b_6t+b_6^2 \big)\lambda_n+2\nu_n-2 t b_6 \right] \Big\} ,
\label{Oparam}
\end{gather}
and
\begin{gather}
  \Theta_{n} = -\frac{b_6\big(1-qb_5^2\big)}{q(1-q)b_5} (x-\lambda_n ) .
\label{Tparam}
\end{gather}
We note that $ \lambda_n$, $\mu_n$, $\nu_n $ satisfy the quadratic relation
\begin{gather}
  \nu_n^2 = (1-q)^2 \lambda_n^2 \mu_n^2 + b_6(b_1 \lambda_n-1) (b_2 \lambda_n-1) (b_3 \lambda_n-1) (b_4 \lambda_n-1)(\lambda_n-t b_6)(b_6 \lambda_n-t) .\!\!\!
\label{Quad}
\end{gather}
\end{proposition}

\begin{proof}
Consistent with the known data, i.e., the degrees, from
(\ref{spectral_coeffEXP:a}), (\ref{spectral_coeffEXP:b}), (\ref{spectral_coeffEXP:c}) we compute the leading
coef\/f\/icients to be
\begin{gather*}
  u_{1} = -\frac{b_6\big(1-qb_5^2\big)}{q(1-q)b_5} ,
\qquad
  v_2 = -\frac{b_6\big(1-b_5^2\big)}{2 (1-q)b_5} ,
\qquad
  w_3 =-\frac{b_6\big(1+b_5^2\big)}{2b_5} ,
\end{gather*}
conf\/irming the relation given by the coef\/f\/icient of $[x^6]$ in (\ref{spectral_Det}), $ w_3^2 = (q-1)^2 v_2^2+b_6^2 $.
In addition we identify the diagonal elements of the $ [x^3] $ coef\/f\/icient of $ A^*_{n} $
\begin{gather*}
 \kappa_{+} \equiv w_3+(q-1)v_2 = -b_5b_6, \qquad
 \kappa_{-} \equiv w_3-(q-1)v_2 = -\frac{b_6}{b_5} .
\end{gather*}
From the coef\/f\/icient of $ [x^0] $ in (\ref{spectral_Det}) we deduce (modulo a sign ambiguity)
\begin{gather*}
 w_0 = b_6t ,
\end{gather*}
and from the coef\/f\/icient of $ [x^1] $ in (\ref{spectral_Det}) we similarly f\/ind
\begin{gather*}
 w_1 = -\tfrac{1}{2} \big[ 1+t(b_1+b_2+b_3+b_4)b_6+b_6^2 \big] .
\end{gather*}
Now utilising the condition $ \nu_{n} = (2W_{n}-W)(\lambda_{n},t) $ we invert this to compute
\begin{gather*}
 w_2 = \frac{1+t (b_1+b_2+b_3+b_4) b_6+b_6^2}{2\lambda_n}+\tfrac{1}{2}b_6\big(b_5+b_5^{-1}\big)\lambda_n
       +\frac{\nu_n-t b_6}{\lambda_n^2} .
\end{gather*}
Proceeding further we infer from the coef\/f\/icient of $ [x^5] $ in (\ref{spectral_Det}) that
\begin{gather*}
 v_1 = \frac{\big(1+b_5^2\big)w_2-\big(b_1^{-1}+b_2^{-1}+b_3^{-1}+b_4^{-1}\big)b_5b_6-b_5\big(1+b_6^2\big)t}{(1-q)\big(1-b_5^2\big)} ,
\end{gather*}
and employing the previous result for $ w_2 $ we derive
\begin{gather*}
  (1-q)\frac{\big(1-b_5^2\big)}{\big(1+b_5^2\big)} v_1 =
  -\frac{b_5b_6\big[ b_1^{-1}+b_2^{-1}+b_3^{-1}+b_4^{-1}+\big(b_6+b_6^{-1}\big)t \big]}{1+b_5^2}
\\
\hphantom{(1-q)\frac{\big(1-b_5^2\big)}{\big(1+b_5^2\big)} v_1 =}{}
+\frac{1+t(b_1+b_2+b_3+b_4)b_6+b_6^2}{2 \lambda_n}
+\frac{b_6\big(1+b_5^2\big)}{2b_5}\lambda_n
+\frac{\nu_n-t b_6}{\lambda_n^2} .
\end{gather*}
This leaves $ v_0 $ to be determined. Imposing the relation $ \mu_{n} =(\Omega_{n}+V)(\lambda_{n},t) $ we
can invert this and f\/ind
\begin{gather*}
 (1-q)\big(1-b_5^2\big)v_0 = (1-q)\big(1-b_5^2\big)\mu_n
   +b_5b_6\big[ b_1^{-1}+b_2^{-1}+b_3^{-1}+b_4^{-1}+\big(b_6+b_6^{-1}\big)t \big]\lambda_n \\
\qquad{}
   -2 b_5b_6\lambda_n^2
  -\tfrac{1}{2}\big(1+b_5^2\big)\big[ 1+t  (b_1+b_2+b_3+b_4 ) b_6+b_6^2 \big]
   +\frac{\big(1+b_5^2\big)}{\lambda_n} (b_6t-\nu_n ) .
\end{gather*}
This concludes our proof.
\end{proof}

\begin{remark}
We observe that the appearance of the quantity $ q^n{b_1b_4b_6} $ with $ n \in \Z_{\geq 0} $ and its replacement by the new parameter $ b_5 $
constitutes a special condition.
This condition is one of the necessary conditions for a member of our particular sequence of classical
solutions to the~$ E^{(1)}_6 $ $q$-Painlev\'e equations, and is built-in by our construction. The
other condition derives from the initial conditions $ n=0 $ in our construction, see \eqref{3termRR} and
following \eqref{spectral_def}.
\end{remark}

From our deformed weight (\ref{DbqJ}) we compute the deformation data polynomials to be
\begin{gather}
 R+\Delta uS  = \frac{1}{b_6} (b_6q t-x ) ,
\label{Ddata1}
\qquad
 R-\Delta uS  =  (q t-b_6x ) .
%\label{Ddata2}
\end{gather}
We can verify that the compatibility relation~(\ref{WV-RS}) is identically satisf\/ied by our spectral and deformation
data polynomials.
We see that this places us in the class $ L=1 $.
We will employ an abbreviation for the dependent variables evaluated at advanced or retarded times, e.g.,
\begin{gather*}
  \lambda_{n}(t) = \lambda_{n}, \qquad \lambda_{n}(qt) = \hat{\lambda}_{n}, \qquad \lambda_{n}\big(q^{-1}t\big) = \check{\lambda}_{n} .
\end{gather*}

In the second stage of our algorithm we parameterise the Cayley transform of the deformation matrix
\begin{gather*}
  B^{*}_{n}
  =  \begin{pmatrix} \mathfrak{R}_{+} & -\mathfrak{P}_{+} \\
                     \mathfrak{P}_{-} &  \mathfrak{R}_{-}
     \end{pmatrix},
  \qquad n \geq 0 ,
\end{gather*}
consistent with known degrees, i.e., $\deg_x \mathfrak{R}_{\pm} = 1$, $\deg_x \mathfrak{P}_{\pm} = 0 $,
so that
\begin{gather*}
 \mathfrak{R}_{\pm} = r_{1,\pm}x+r_{0,\pm} ,
\qquad
 \mathfrak{P}_{\pm} = p_{\pm} .
\end{gather*}

\begin{lemma}
Let us assume $ b_6 \neq 0 $ and $ b_5 \neq q^{-1/2}, q^{1/2} $.
Then the off-diagonal components of the deformation matrix are given by
\begin{gather}
   p_{+} = -\hat{a}_nr_{1,-}+a_nr_{1,+} ,
\label{pr:a}
\\
   p_{-} = -a_nr_{1,-}+\hat{a}_nr_{1,+} .
\label{pr:b}
\end{gather}
\end{lemma}
\begin{proof}
We resolve the $ A $-$ B $ compatibility relation (\ref{A-B}) into monomials of $ x $. Examining the~$ x^7 $ coef\/f\/icient of the
$ (1,2) $ and $ (2,1) $ components yields (\ref{pr:a}) and (\ref{pr:b}) respectively.
\end{proof}

\begin{lemma}
Let us assume $ b_5 \neq q^{1/2} $, $ a_{n}, \hat{a}_{n} \neq 0 $ and $ \lambda_{n} \neq b_6 t, b_6^{-1} t $.
Then the spectral and deformation matrices satisfy the following residue formulae
\begin{gather}
 \mathfrak{R}_{-}(b_6 qt,t)+\frac{\mathfrak{W}_{+}(b_6 qt,qt)}{\mathfrak{T}_{+}(b_6 qt,qt)} \mathfrak{P}_{+}(b_6 qt,t) = 0 ,
\label{EQ1} \\
 \mathfrak{R}_{-}(b_6^{-1}qt,t)+\frac{\mathfrak{W}_{+}(b_6^{-1}qt,qt)}{\mathfrak{T}_{+}(b_6^{-1}qt,qt)} \mathfrak{P}_{+}\big(b_6^{-1}qt,t\big) = 0 ,
\label{EQ2}
\end{gather}
and
\begin{gather}
 \mathfrak{R}_{+}(b_6 qt,t)+\frac{\mathfrak{W}_{-}(b_6 t,t)}{\mathfrak{T}_{+}(b_6 t,t)} \mathfrak{P}_{+}(b_6 qt,t) = 0 ,
\label{EQ3}\\
 \mathfrak{R}_{+}(b_6^{-1}qt,t)+\frac{\mathfrak{W}_{-}(b_6^{-1}t,t)}{\mathfrak{T}_{+}(b_6^{-1}t,t)} \mathfrak{P}_{+}\big(b_6^{-1}qt,t\big) = 0 .
\label{EQ4}
\end{gather}
\end{lemma}

\begin{proof}
In this step we compute the residues of the $ A $-$ B $ compatibility relation, with respect to~$ x $,
at the zeros and poles of
\begin{gather}
 \chi(x,t)=
 \frac{\left(x-qb_6 t\right) \left(b_6 x-q t\right)}{q \left(x-b_6 t\right) \left(b_6 x-t\right)} .
\end{gather}
From the residue of (\ref{S+D:a}) at the zero $ x=b_6 qt $ we deduce~(\ref{EQ1}),
and from the same equation at the zero $ x=b_6^{-1} qt $ we deduce~(\ref{EQ2}).
From the residue of~(\ref{S+D:d}) at the pole $ x=b_6 t $ we deduce~(\ref{EQ3}), and
from the same equation at the pole $ x=b_6^{-1} t $ we deduce~(\ref{EQ4}).
\end{proof}

\begin{remark}
Although the above proof appealed to the vanishing of the right-hand side of one of the compatibility
conditions, namely (\ref{S+D:a}), at either of the two zeros of $ \chi $, in fact under these conditions
the right-hand sides of all the other compatibility conditions, i.e.~(\ref{S+D:b}),~(\ref{S+D:c}), and~(\ref{S+D:d}),
also vanish. This is because $ \chi=0 $ implies $ \big(R^2-\Delta u^2S^2\big)\big(b_6^{\pm 1}qt;t\big)=0 $ and
$ \big(W^2-\Delta y^2V^2\big)\big(b_6^{\pm 1}qt;qt\big)=0 $, and furthermore the spectral and deformation matrices satisfy
the determinantal identities
\begin{gather*}
  \det A^{*}_n = \mathfrak{W}_{+}\mathfrak{W}_{-}+\mathfrak{T}_{+}\mathfrak{T}_{-} = W^2-\Delta y^2V^2,
\\
  \det B^{*}_n = \mathfrak{R}_{+}\mathfrak{R}_{-}+\mathfrak{P}_{+}\mathfrak{P}_{-} = \frac{a_n}{\hat{a}_n}\big( R^2-\Delta u^2S^2 \big).
\end{gather*}
Therefore under the specialisations $ x=b_6^{\pm 1}qt $ the right-hand sides of (\ref{S+D:b}), (\ref{S+D:c}), and (\ref{S+D:d})
are proportional to the right-hand side of~(\ref{S+D:a}), and the vanishing of the latter implies
the vanishing of the former. In this way we ensure that all components of the $A$-$B$ compatibility
vanish under the single condition. A similar observation applies to the left-hand sides of the
compatibility relations at the zeros of $ \chi^{-1} $, i.e.~$ x=b_6^{\pm 1}t $.
\end{remark}

We introduce our f\/irst change of variables, $ \mu_n, \nu_n \mapsto \mathit{z}_{\pm} $, via the relations
\begin{gather}
   \nu_n = \frac{1}{2}\lambda_{n} [\kappa _{+}\mathit{z}_{+}+\kappa _{-}\mathit{z}_{-} ] ,
\label{MNsubs:a}
\\
   \mu_n = \frac{1}{2(q-1)} [\kappa _{+}\mathit{z}_{+}-\kappa _{-}\mathit{z}_{-} ] .
\label{MNsubs:b}
\end{gather}
The new variables satisfy an identity corresponding to \eqref{Quad} which reads
\begin{gather*}
\kappa _{+}\kappa _{-}\mathit{z}_{+}\mathit{z}_{-}
= \frac{1}{\lambda_{n}^2}\big[ W^2-\Delta y^2V^2 \big](\lambda_{n},t) .
\end{gather*}

Next we subtract (\ref{EQ2}) from (\ref{EQ1}), in order to eliminate both $ \mathit{z}_{-} $ and $ r_{0,-} $.
This yields
\begin{gather}
 \frac{qp_{+}}{\big(1-q b_5^2\big)\hat{a}_n}\left[ -\frac{b_5\big(b_5 qt\hat{\lambda}_n-1\big)}{b_6\hat{\lambda}_n}
            +\frac{qb_5^2 t}{\big(\hat{\lambda}_n-b_6q t\big) \big(b_6\hat{\lambda}_n-q t\big) }\hat{\mathit{z}}_{+} \right]
		  +q t\frac{1}{b_6} r_{1,-} = 0 .
\label{EQ5}
\end{gather}
This result motivates the def\/inition of the new variable $ \mathpzc{Z} $
\begin{gather*}
  \mathpzc{Z} = \left.-\frac{b_5 b_6qt}{(\hat{\lambda}_n-b_6 qt)(b_6 \hat{\lambda}_n-qt)}\hat{\mathit{z}}_{+}
   +\frac{b_5 qt\hat{\lambda}_n-1}{\hat{\lambda}_n}\right|_{t\to  q^{-1}t} .
\end{gather*}
\begin{definition}
In terms of this new variable $ \mathpzc{Z} $ we have
\begin{gather}
   \mathit{z}_{+} = \frac{1}{b_5 b_6 t}\frac{(\lambda_n-t b_6)(b_6 \lambda_n-t) [(b_5 t -\mathpzc{Z})\lambda_n-1 ]}{\lambda_n} ,
\label{Zsubs:a}
\\
   \mathit{z}_{-} = b_5 t \frac{(b_1 \lambda_n-1)(b_2 \lambda_n-1)(b_3 \lambda_n-1)(b_4 \lambda_n-1)}
                               {\lambda_n  [(b_5 t-\mathpzc{Z})\lambda_n-1 ]} .
\label{Zsubs:b}
\end{gather}
\end{definition}
Our f\/inal rewrite of the dependent variables is
\begin{gather}
 \lambda_n(t) \to g(t) ,
\label{FGsubs:a}
\\
 \mathpzc{Z}(t) \to  b_5 t -f(q^{-1}t) .
\label{FGsubs:b}
\end{gather}

We are now in a position to undertake the third stage of our derivation.
The f\/irst of the evolution equations is given in the following result.
\begin{proposition}
Let us assume that $ q \neq 0 $, $ b_5 \neq q^{-1/2} $, $ t \neq 0 $, $ g \neq 0, b_6t, b_6^{-1}t $ and $ a_{n} \neq 0 $.
The variables~$ f$,~$g $ satisfy the first-order $q$-difference equation
\begin{gather}
\big(g\check{f}-1\big)(gf-1)
 = t^2\frac{\big(g-b_1^{-1}\big)\big(g-b_2^{-1}\big) \big(g-b_3^{-1}\big)\big(g-b_4^{-1}\big)}
           {(g-b_6t)\big(g-b_6^{-1}t\big)} .
\label{1stEvol}
\end{gather}
This evolution equation is identical to the second equation of equation~{\rm (4.15)} of Kajiwara et al.~{\rm \cite{KMNOY_2004}} and to the second equation
of equation~{\rm (3.23)} of Kajiwara et al.~{\rm \cite{KMNOY_2005}}.
\end{proposition}

\begin{proof}
Subtract (\ref{EQ4}) from (\ref{EQ3}) in order to eliminate both $ \mathit{z}_{+} $ and $ r_{0,+} $.
This yields the relation
\begin{gather}
 \frac{q p_{+}}{\big({-}1+q b_5^2\big) a_n}\left[ \frac{t}{(b_6t-\lambda_n) (-t+b_6 \lambda_n)}\mathit{z}_{-}
                                                                     -\frac{(b_5-t \lambda_n)}{b_6 \lambda_n} \right]
  +q t \frac{1}{b_6}r_{1,+} = 0 .
\label{EQ6}
\end{gather}
Thus we have two dif\/ferent ways of computing the ratio of $ r_{1,+} $ to $ r_{1,-} $; on the one hand we have from (\ref{EQ5})
\begin{gather}
 r_{1,+} = -\frac{[ b_5 (q b_5 t -\hat{\mathpzc{Z}})-t ] \hat{a}_n}{b_5 \hat{\mathpzc{Z}} a_n} r_{1,-} ,
\label{aux1}
\end{gather}
whereas using (\ref{EQ6}) we have
\begin{gather*}
 \frac{b_5 a_n}{\hat{a}_n}\frac{r_{1,+}}{r_{1,-}} =
    \Big[ b_5 b_6  t^2 \left(b_1 \lambda_n-1\right) \left(b_2 \lambda_n-1\right) \left(b_3 \lambda_n-1\right) \left(b_4 \lambda_n-1\right)
\nonumber\\
  \hphantom{\frac{b_5 a_n}{\hat{a}_n}\frac{r_{1,+}}{r_{1,-}} = }\qquad{}
   +\left(b_5-t \lambda_n\right)\left(\lambda_n-t b_6\right) \left(b_6 \lambda_n-t\right)\left[\lambda_n\left(b_5 t-\mathpzc{Z}\right)-1\right] \Big]
\nonumber\\
\hphantom{\frac{b_5 a_n}{\hat{a}_n}\frac{r_{1,+}}{r_{1,-}} = }
 \div \Big[ t^2 b_6 \left(b_1 \lambda_n-1\right) \left(b_2 \lambda_n-1\right) \left(b_3 \lambda_n-1\right) \left(b_4 \lambda_n-1\right)
\nonumber\\
\hphantom{\frac{b_5 a_n}{\hat{a}_n}\frac{r_{1,+}}{r_{1,-}} = }\qquad{}
 -\left(q b_5 t \lambda_n-1\right)\left(\lambda_n-t b_6\right) \left(b_6 \lambda_n-t\right)\left[\lambda_n\left(b_5  t-\mathpzc{Z}\right)-1\right] \Big] .
\end{gather*}
Equating these two forms gives (\ref{1stEvol}).
\end{proof}

The second evolution equation, to be paired with the f\/irst~(\ref{1stEvol}) as a coupled system, is given next.
\begin{proposition}
Let us make the following assumptions: $ t \neq 0 $, $ b_5\neq 1, q^{-1/2}, q^{-1} $, $ f \neq 0 $, $ b_5f \neq t $,
$ g \neq 0, \hat{g} \neq 0 $ and $ a_{n} \neq 0 $.
In addition let us assume that the condition
\begin{gather*}
 \hat{g}\neq \frac{1-q b_5  t  g-q b_5^2+q b_5^2 f g}{f-b_5 qt} ,
\end{gather*}
holds. The variables $ f$, $g $ satisfy the first-order $q$-difference equation
\begin{gather}
(f \hat{g}-1)(f g-1) =
 qt^2\frac{(f-b_1) (f-b_2)(f-b_3) (f-b_4)}
                           {(f-b_5qt)\big(f-b_5^{-1}t\big)} .
\label{2ndEvol}
\end{gather}
This evolution equation is the same as the first equation of equation~{\rm (4.15)} in Kajiwara et al.~{\rm \cite{KMNOY_2004}}
and the first equation of equation~{\rm (3.23)} in Kajiwara et al.~{\rm \cite{KMNOY_2005}}, both subject to typographical corrections.
\end{proposition}
\begin{proof}
Cross multiplying the relations (\ref{EQ1}), (\ref{EQ2}), (\ref{EQ3}), (\ref{EQ4}) we can eliminate all refe\-ren\-ce
to the deformation matrix and deduce the identity
\begin{gather}
 \frac{\mathfrak{W}_{+}(b_6q t,q t)}{\mathfrak{T}_{+}(b_6q t,q t)}
 \frac{\mathfrak{W}_{-}(b_6t,t)}{\mathfrak{T}_{+}(b_6t,t)}
 = \frac{\mathfrak{W}_{+}\big(q b_6^{-1}t,q t\big)}{\mathfrak{T}_{+}\big(q b_6^{-1}t,q t\big)}
   \frac{\mathfrak{W}_{-}\big(b_6^{-1}t,t\big)}{\mathfrak{T}_{+}\big(b_6^{-1}t,t\big)} .
\label{Idprod}
\end{gather}
Into this identity we employ the following evaluations for the advanced and retarded values of~$ \mathit{z}_{\pm} $
\begin{gather*}
 \hat{\mathit{z}}_{+} = \frac{q^{-1}}{b_5 b_6 t} \frac{(f \hat{g}-1) (\hat{g}- b_6 qt) (\hat{g} b_6-q t)}{\hat{g}} ,
\\
 \hat{\mathit{z}}_{-} = q b_5 t \frac{(\hat{g} b_1-1) (\hat{g} b_2-1) (\hat{g} b_3-1)(\hat{g} b_4-1)}{\hat{g} (f \hat{g}-1)} ,
\\
 \mathit{z}_{+} = \frac{t}{ b_5}\frac{(g b_1-1) (g b_2-1)(g b_3-1) (g b_4-1t)}{g (f g-1)} ,
\\
 \mathit{z}_{-} = \frac{b_5}{b_6 t}\frac{(f g-1)(g-b_6 t) (g b_6-t)}{g} .
\end{gather*}
We f\/ind that this relation factorises into two non-trivial factors, the f\/irst of which is proportional to
\begin{gather*}
 \hat{g} - \frac{1-q b_5 t  g-q b_5^2 +q b_5^2 f g}{f-q b_5 t} .
\end{gather*}
Assuming this is non-zero our evolution equation is then the remaining factor of~(\ref{Idprod})
\begin{gather*}
(f \hat{g}-1)(f g-1) =
 q b_5 t^2 \frac{(f-b_1) (f-b_2)(f-b_3) (f-b_4)}
                           {(f-q b_5 t ) (b_5 f -t)} ,
\end{gather*}
or alternatively (\ref{2ndEvol}).
\end{proof}

Lastly we have an auxiliary evolution equation which controls the normalisation of the orthogonal polynomial system.
\begin{proposition}
Let us assume $ b_6 \neq 0 $, $ b_5 \neq q^{-1/2} $, $ f \neq b_5qt, b_5^{-1}t $ and $ \gamma_n \neq 0 $ for $ n \geq 0 $.
The leading coefficient of the polynomials or second-kind solutions $($see \eqref{ops_pExp}, \eqref{ops_eExp}$)$ satisfy the
first-order $q$-difference equation
\begin{gather}
  \left(\frac{\hat{\gamma}_n}{b_6\gamma _n}\right)^2 = \frac{f-b_5^{-1}t}{f-b_5qt} .
\label{3rdEvol}
\end{gather}
\end{proposition}
\begin{proof}
Using the leading order, i.e., the $ [x] $ terms, in the expansions (\ref{deform_coeff_Exp:a}), (\ref{deform_coeff_Exp:d})
with def\/initions (\ref{deform_def:abc}) we can compute~$ r_{1,+} $. However by using these
same expansions to compute~$ r_{1,-} $ and the equation~(\ref{aux1}), which relates these two quantities, we have
an alternative expression for~$ r_{1,+} $. Equating these expressions then gives~(\ref{3rdEvol}).
\end{proof}

We conclude our discussion by summarising our results for the spectral and deformation matrices
in terms of the $ f_{n}$, $g_{n} $ variables. Henceforth we will restore the index $ n $ on all our variables.
The form of the spectral matrix is given in the following proposition.
\begin{proposition}
Assume that $ |q| \neq 1 $, $ b_6 \neq 0 $, $ b_5^2 \neq q^{-1},1,q $ and $ a_{n} \neq 0 $.
The spectral matrix elements \eqref{spectral_prod}, \eqref{spectral_def}, \eqref{AC_defn} are given by
\begin{gather*}
  \mathfrak{T}_{n,+} = \frac{b_6}{b_5}\big(q^{-1}-b_5^2\big) a_n x (x-g_{n} ) ,
\qquad
  \mathfrak{T}_{n,-} = \frac{b_6}{b_5}\big(1-q^{-1} b_5^2\big) a_n x (x-g_{n-1} ) ,
\end{gather*}
and
\begin{gather}
\frac{\mathfrak{W}_{n,+}}{x(x-g_{n})}
  =   -x b_5 b_6+\frac{b_6}{\big(1-b_5^2\big)}\left[ -\frac{b_5^2}{t}+b_5\left(\frac{1}{b_1}+\frac{1}{b_2}+\frac{1}{b_3}+\frac{1}{b_4}\right) \right]
\nonumber\\
  \hphantom{\frac{\mathfrak{W}_{n,+}}{x(x-g_{n})}=}{}
      -\frac{ b_6(b_5f_{n}-t)}{(1-b_5^2)t}\left[ \frac{g_{n}(t-f_{n} b_5)}{f_{n}}+\frac{t b_5 \big(1+b_6^2\big)}{b_6} \right]
\nonumber\\
\hphantom{\frac{\mathfrak{W}_{n,+}}{x(x-g_{n})}=}{}
 +\frac{b_6t}{(1-b_5^2)} \left[ \frac{b_5^2}{g_{n}^2}-\frac{1-b_5^2}{x g_{n}}-(b_1+b_2+b_3+b_4)b_5^2\frac{1}{g_{n}}+ b_5^2\frac{f_{n}}{g_{n}}-\frac{g_{n}}{f_{n}} \right.
\nonumber\\
 \left.
\hphantom{\frac{\mathfrak{W}_{n,+}}{x(x-g_{n})}=}{}
       +\frac{ (1-g_{n} b_1) (1-g_{n} b_2) (1-g_{n} b_3) (1-g_{n} b_4) \big(g_{n}-x b_5^2\big)}{(1-f_{n} g_{n}) g_{n}^2 (x-g_{n})} \right] ,
\label{eq:Anstar1}
\end{gather}
or
\begin{gather*}
\frac{\mathfrak{W}_{n,+}}{x(x-g_{n})}
  =   -\frac{b_6t}{\big(1-b_5^2\big)}\frac{(f_{n}-b_1)(f_{n}-b_2)(f_{n}-b_3)(f_{n}-b_4)\big(1-b_5^2f_{n}x\big)}{f_{n}^2(1-f_{n}g_{n})(1-f_{n}x)}
\nonumber\\
\hphantom{\frac{\mathfrak{W}_{n,+}}{x(x-g_{n})}=}{}
  +\frac{b_6 t(1-x b_1)(1-x b_2)(1-x b_3)(1-x b_4)}{x(1-xf_{n})(x-g_{n})} -\frac{b_6(b_5f_{n}-t)}{f_{n}}x
\nonumber\\
\hphantom{\frac{\mathfrak{W}_{n,+}}{x(x-g_{n})}=}{}
    -\frac{b_6(b_5f_{n}-t)}{b_5\big(1-b_5^2\big)t}\left\{ \frac{b_5^2\big(1+b_6^2\big)t}{b_6}+\frac{b_5g_{n}}{f_{n}}(t-b_5f_{n}) \right.
\nonumber\\
 \left.
\hphantom{\frac{\mathfrak{W}_{n,+}}{x(x-g_{n})}=}{}
    +\frac{b_5}{f_{n}^2} \left[ t+b_5f_{n}-t f_{n} \left(\frac{1}{b_1}+\frac{1}{b_2}+\frac{1}{b_3}+\frac{1}{b_4}\right) \right]
                                    \right\} ,
\end{gather*}
and
\begin{gather}
\frac{t\mathfrak{W}_{n,-}}{x(x-g_{n})}
     =   \frac{(1-xf_{n})(x-b_6t)(x b_6-t)}{x(x-g_{n})}-\frac{b_6}{1-b_5^2}\frac{(t-b_1 b_5)(t-b_2 b_5)(t-b_3 b_5)(t-b_4 b_5)}{b_5(t g_{n}-b_5)}
\nonumber\\
\hphantom{\frac{t\mathfrak{W}_{n,-}}{x(x-g_{n})}=}{}
 +\frac{b_6(b_5f_{n}-t)}{b_5\big(1-b_5^2\big)}\left[ \big(1-b_5^2\big)x-\frac{t^2}{g_{n}}-g_{n} b_5^2+\frac{b_5^2\big(1+b_6^2\big)t}{b_6} \right.
\nonumber\\
          \left.
\hphantom{\frac{t\mathfrak{W}_{n,-}}{x(x-g_{n})}=}{}
           -\frac{b_5t^2 (1-g_{n} b_1)(1-g_{n} b_2)(1-g_{n} b_3)(1-g_{n} b_4)}{g_{n}(1-f_{n}g_{n})(tg_{n}-b_5)} \right] ,
\label{eq:Anstar2}
\end{gather}
or
\begin{gather*}
\frac{t\mathfrak{W}_{n,-}}{x(x-g_{n})}
    =   \frac{(1-x f_{n})(x-b_6 t)(b_6 x-t)}{x(x-g_{n})}+\frac{b_6t^2}{\big(1-b_5^2\big)}\frac{(f_{n}-b_1)(f_{n}-b_2)(f_{n}-b_3)(f_{n}-b_4)}{f_{n}^2(1-f_{n}g_{n})}
\nonumber\\
\hphantom{\frac{t\mathfrak{W}_{n,-}}{x(x-g_{n})}=}{}
 + \frac{b_6(b_5f_{n}-t)}{b_5\big(1-b_5^2\big)}\left\{ \big(1-b_5^2\big)x+\frac{b_5^2\big(1+b_6^2\big)t}{b_6}+\frac{b_5g_{n}}{f_{n}}(t-b_5 f_{n}) \right.
\nonumber\\
\left.
\hphantom{\frac{t\mathfrak{W}_{n,-}}{x(x-g_{n})}=}{}
 +\frac{b_5}{f_{n}^2}\left[ t+b_5f_{n}-tf_{n}\left(\frac{1}{b_1}+\frac{1}{b_2}+\frac{1}{b_3}+\frac{1}{b_4}\right) \right] \right\} .
\end{gather*}
\end{proposition}

\begin{proof}
This follows from applying the transformations (\ref{MNsubs:a}), (\ref{MNsubs:b}), (\ref{Zsubs:a}), (\ref{Zsubs:b}) and
(\ref{FGsubs:a}), (\ref{FGsubs:b}) successively to~(\ref{Tparam}), (\ref{Wparam}) and~(\ref{Oparam}).
The alternative forms arise from applying partial fraction expansions with respect to either of
$ f_{n} $ or $ g_{n} $.
\end{proof}

The deformation matrix is summarised in the next result.
\begin{proposition}
Assume that $ |q| \neq 1 $, $ b_6 \neq 0 $, $ b_5^2 \neq q^{-1},1,q $.
The deformation matrix elements are given by
\begin{gather}
 \mathfrak{R}_{n,+} = \frac{\hat{\gamma}_{n}}{b_6\gamma_{n}} \left\{ x +\frac{1}{1-b_5^2}\left[
   -\frac{q \big(t^2 b_6+g_{n}^2 b_5^2 b_6-t g_{n} b_5^2\big(1+b_6^2\big)\big)}{g_{n}b_6}\right.\right.
\nonumber\\
\hphantom{\mathfrak{R}_{n,+} =}{}
   +\frac{q t^2 b_5(1-b_1g_{n})(1-b_2g_{n})(1-b_3g_{n})(1-b_4g_{n})}{g_{n} (f_{n} g_{n}-1) (t g_{n}-b_5)}
\nonumber\\
\left.\left.
\hphantom{\mathfrak{R}_{n,+} =}{}   -\frac{q(t-b_1b_5)(t-b_2b_5)(t-b_3b_5)(t-b_4b_5)}{(t g_{n}-b_5) (f_{n} b_5-t)} \right] \right\} ,
\label{DMatrix11}
\end{gather}
and
\begin{gather*}
 \mathfrak{R}_{n,-} = \frac{b_6\gamma_{n-1}}{\hat{\gamma}_{n-1}} \left\{ x +\frac{1}{1-b_5^2}\left[
 \frac{q\big(t^2 b_6+g_{n}^2 b_5^2 b_6-t g_{n}\big(1+b_6^2\big)\big)}{g_{n}b_6}\right.\right.
\nonumber\\
\hphantom{\mathfrak{R}_{n,-} =}{}
  -\frac{q t^2 b_5(1-b_1g_{n})(1-b_2g_{n})(1-b_3g_{n})(1-b_4g_{n})}{g_{n} (f_{n} g_{n}-1) (t g_{n}-b_5)}
\nonumber\\
\left.\left.
\hphantom{\mathfrak{R}_{n,-} =}{}
 +\frac{q(t-b_1b_5)(t-b_2b_5)(t-b_3b_5)(t-b_4b_5)}{(t g_{n}-b_5) (f_{n} b_5-t)} \right] \right\} .
%\label{DMatrix22}
\end{gather*}
Furthermore
\begin{gather*}
 \mathfrak{P}_{n,+} = a_{n}\left[ \frac{\hat{\gamma}_{n}}{b_6\gamma_{n}}-\frac{b_6\gamma_{n}}{\hat{\gamma}_{n}} \right] ,
\qquad
 \mathfrak{P}_{n,-} = a_{n}\left[ \frac{\hat{\gamma}_{n-1}}{b_6\gamma_{n-1}}-\frac{b_6\gamma_{n-1}}{\hat{\gamma}_{n-1}} \right] .
\end{gather*}
\end{proposition}

\begin{proof}
Using the leading orders in the expansions (\ref{deform_coeff_Exp:a}), (\ref{deform_coeff_Exp:d}), i.e., the $ [x] $ terms,
with def\/initions (\ref{deform_def:abc}) we deduce
\begin{gather*}
 r_{1,+} = \frac{\hat{\gamma}_{n}}{b_6\gamma_{n}} ,
\qquad
 r_{1,-} = \frac{b_6\gamma_{n-1}}{\hat{\gamma}_{n-1}} .
\end{gather*}
Using the leading orders in the expansions (\ref{deform_coeff_Exp:b}), i.e., the $ [x^0] $ terms,
with def\/inition~(\ref{deform_def:abc}) and~(\ref{Defm_offD}) we deduce
\begin{gather*}
 p_{+} = a_{n}\left[ \frac{\hat{\gamma}_{n}}{b_6\gamma_{n}}-\frac{b_6\gamma_{n}}{\hat{\gamma}_{n}} \right] ,
\qquad
 p_{-} = a_{n}\left[ \frac{\hat{\gamma}_{n-1}}{b_6\gamma_{n-1}}-\frac{b_6\gamma_{n-1}}{\hat{\gamma}_{n-1}} \right] .
\end{gather*}
Using the coef\/f\/icient of the $ [x^7] $ term in the $ (1,1) $ element of the $ A $-$ B $ compatibility rela\-tions~(\ref{spectral+deform:b}), along with the solution of~(\ref{1stEvol}) for $ f_{n}(q^{-1}t) $ and~(\ref{2ndEvol})
for~$ g_{n}(qt) $ we deduce
\begin{gather*}
  \big(1-b_5^2\big)\frac{r_{0,+}}{r_{1,+}}
  = q b_5\left(\frac{t}{f_{n}}-b_5\right)g_{n}+\left(\frac{b_5 qt}{f_{n}}-1\right)\hat{g}_{n}\\
  \hphantom{\big(1-b_5^2\big)\frac{r_{0,+}}{r_{1,+}}=}{}
          +\big[ 1+q b_5^2-b_5\big(b_1^{-1}+b_2^{-1}+b_3^{-1}+b_4^{-1}\big)qt \big]\frac{1}{f_{n}}+\frac{\big(1+b_6^2\big)}{b_6}b_5^2 qt ,
\end{gather*}
whereas if we examine the $ [x^7] $ term in the $ (2,2) $ element of the $ A $-$ B $ compatibility relations in the
same way then we f\/ind
\begin{gather*}
 \big(1-b_5^2\big)\frac{r_{0,-}}{r_{1,-}}
  = q b_5\left(b_5-\frac{t}{f_{n}}\right)g_{n}+\left(1-\frac{b_5 qt}{f_{n}}\right)\hat{g}_{n}\\
\hphantom{\big(1-b_5^2\big)\frac{r_{0,-}}{r_{1,-}}=}{}
          +\big[ b_5\big(b_1^{-1}+b_2^{-1}+b_3^{-1}+b_4^{-1}\big)qt-q b_5^2-1 \big]\frac{1}{f_{n}}-\frac{\big(1+b_6^2\big)}{b_6}qt  .
\end{gather*}
Into both of these expressions we can employ~(\ref{2ndEvol}) for~$ \hat{g}_n $ and make a partial
fraction expansion with respect to~$ f_{n} $.
\end{proof}

\section{Reconciliation with the Lax pairs of Sakai and Yamada}\label{Reconcile}

\subsection{Sakai Lax Pair}

In \cite{Sa_2006} Sakai constructed a Lax pair for the $E^{(1)}_6$ $q$-Painlev\'e equations using a degeneration of
a two-variable case of the Garnier system based upon the Lax pairs for the $D^{(1)}_5$ $q$-Painlev\'e system~\cite{Sa_2005}.
Subsequently Murata~\cite{Mu_2009} gave more details for this Lax pair. We intend to establish a correspondence between
our Lax pair and that of Sakai. We will carry this out in a sequence of simple steps rather than as a single step as this will
reveal how similar they are.

Our f\/irst step is to give a variation on the parameterisation of the spectral and deformation matrices to that
given in Section~\ref{Big_qJacobi}. In this alternative formulation, we seek a spectral matrix $\tilde{A}(x;t)$ (actually
identical to the Cayley transform $A_n^*$) with the specif\/ications
\begin{gather}
\tilde{A}(x;t) = A_{0} + A_{1}x + A_{2}x^2 + A_{3}x^3,	
\label{eq:new1}
\end{gather}
and
\begin{enumerate}\itemsep=0pt
\item the determinant is
\[
b_6 \left(b_1 x-1\right) \left(b_2 x-1\right) \left(b_3 x-1\right) \left(b_4 x-1\right) \left(x-b_6 t\right) \left(b_6 x-t\right) ,
\]

\item $A_{3}$ is diagonal with entries $\kappa_1= -b_5b_6$ and $\kappa_2 = -b_6/b_5$,
\item $A_{0} = b_6 t \mathbbmss{1}$,
\item the root of the (1,2) entry of $\tilde{A}(x;t)$ with respect to $x$ is $\lambda$,
\item $\tilde{A}(\lambda;t)$ is lower triangular with diagonal entries $-b_5b_6\lambda z_{+}$ and $-\frac{b_6\lambda z_{-}}{b_5}$ where $b_6^2z_-z_+\lambda^2 = \det \tilde{A}(\lambda;t)$.
\end{enumerate}
Any such matrix is in the general form
\begin{gather*}
\tilde{A}(x;t) = tb_6I - \begin{pmatrix}  b_5 b_6 x  [ z_1+(x-\alpha )  (x-\lambda  ) ] & \dfrac{b_6 w x \left(x-\lambda \right)}{b_5} \\
 \dfrac{b_5 b_6 x (x \gamma +\delta ) }{w} & \dfrac{b_6 x [ z_2+(x-\beta ) (x-\lambda ) ] }{b_5}
                           \end{pmatrix},
%\label{An_new}
\end{gather*}
where the properties specify the variables
\begin{gather}
\big(1 - b_5^2\big) \alpha = \frac{1}{b_1} + \frac{1}{b_2} + \frac{1}{b_3} + \frac{1}{b_4} + \left( \frac{1}{b_6} + b_6 \right)t
\nonumber\\
\hphantom{\big(1 - b_5^2\big) \alpha =}{}
	- (b_1 + b_2 + b_3 + b_4)b_5 \frac{t}{\lambda} -b_5\left( \frac{1}{b_6} + b_6 \right) \frac{1}{\lambda}
		 + \frac{b_5^2 z_1}{\lambda} +  \frac{z_2}{\lambda} - 2 \lambda,
\label{alpha_new}
\\
\big(1 - b_5^2\big) \beta = - \left(  \frac{1}{b_1} + \frac{1}{b_2} + \frac{1}{b_3} + \frac{1}{b_4}  \right) b_5^2
	- \left( \frac{1}{b_6} + b_6 \right)b_5^2 t
\nonumber\\
\hphantom{\big(1 - b_5^2\big) \beta =}{}
      + (b_1 + b_2 + b_3 + b_4)b_5 \frac{t}{\lambda} + b_5\left( \frac{1}{b_6} + b_6 \right) \frac{1}{\lambda}
			- \frac{b_5^2 z_1}{\lambda} -  \frac{z_2}{\lambda} + 2 b_5^2 \lambda,
\label{beta_new}
\\
\gamma = -  ( b_3 b_4 + b_2 b_3 + b_2 b_4 + b_1 b_2 + b_1 b_3  + b_1 b_4  )
	-  \left(  \frac{1}{b_1} + \frac{1}{b_2} + \frac{1}{b_3} + \frac{1}{b_4}  \right) \left( \frac{1}{b_6} + b_6 \right) t \nonumber\\
\hphantom{\gamma =}{}
		-t^2+ \alpha \beta +z_1 + z_2 + 2 (\alpha + \beta) \lambda + \lambda^2 ,
\label{gamma_new}
\\
\delta = b_1 + b_2 + b_3 + b_4 +\left[   ( b_3 b_4 + b_2 b_3 + b_2 b_4 + b_1 b_2 + b_1 b_3  + b_1 b_4  )
	\left(\! \frac{1}{b_6} + b_6\! \right) - \left(\! \frac{1}{b_5} + b_5 \!\right) \right]t\nonumber \\
\hphantom{\delta =}{}
+ \left( \frac{1}{b_1} + \frac{1}{b_2} + \frac{1}{b_3} + \frac{1}{b_4} \right) t^2 - z_1(\beta + \lambda)
		-z_2 (\alpha + \lambda)  + (-2\alpha \beta + \gamma) \lambda - (\alpha + \beta) \lambda^2 .\!\!\!
\label{delta_new}
\end{gather}
The $z_1$ and $z_2$ are related to $z_{\pm}$ by
\begin{gather}
    z_1  = z_+ + \frac{t}{b_5 \lambda}  \qquad  \mbox{and} \qquad z_2 = z_{-} + \frac{ tb_5}{\lambda} .
\label{z_new}
\end{gather}
In addition
\begin{gather*}
    w = \frac{1 - q b_5^2}{q} a_n .
\end{gather*}
We seek a deformation matrix $\tilde{B}(x;t)$ of the form
\begin{gather}
   \tilde{B}(x;t) = \frac{x}{(x - b_6 qt) \big( x - b_6^{-1} qt\big)} (x \mathbbmss{1} + B_{0}), 	
\qquad\mbox{where}\qquad
          B_{0} = \begin{bmatrix} r_{1,1} & r_{1,2} \\
				    r_{2,1} & r_{2,2}
		  \end{bmatrix} .
\label{eq:Sakai2}
\end{gather}
This leads to the compatibility relation
\begin{gather}
  \tilde{B}(qx;t) \tilde{A}(x;t) = \tilde{A}(x;qt) \tilde{B}(x;t).
\label{eq:new3}
\end{gather}
This relation is just a rewriting of \eqref{spectral+deform:b} whereby all the factors of $\chi$ are placed into
the denominator of $\tilde{B}$ by the above def\/inition.

\begin{lemma} \label{lem:new1}
The overdetermined system \eqref{eq:new3}, with \eqref{eq:new1} and \eqref{eq:Sakai2} is satisfied if the coupled
$E^{(1)}_6$ $q$-Painlev\'e equations \eqref{1stEvol} and \eqref{2ndEvol} are satisfied.
\end{lemma}

\begin{proof}
Examining the coef\/f\/icient of $x^6$ in the numerator of the (1,2) entry of \eqref{eq:new3} we f\/ind
\begin{gather}
r_{1,2} = \frac{q}{1 -qb_5^2}  ( \hat{w} - w  ) .
\label{eq:new6}
\end{gather}
Now we seek two alternative expressions for $r_{1,2}$~-- one involving quantities at the advanced time $qt$ and
another involving those at the unshifted time $t$. The f\/irst of these is found from solving for the (1,2) entry
of the residue of~\eqref{eq:new3} at $x = b_6 q t$ simultaneously with the (1,2) entry of the residue of~\eqref{eq:new3}
at $x = b_6^{-1} q t$. This yields
\begin{gather}
r_{1,2} = \frac{ - q t \hat{w} \big(q t b_6 - \hat{\lambda}\big) \big( q t -b_6 \hat{\lambda} \big)}
		{b_5 \left\{ q t \big(1 - b_5 b_6 \hat{z}_1\big) + \big(q t b_6 - \hat{\lambda}\big)  \big[ b_6 + q t b_5 \big( q t - b_6 \hat{\lambda} \big) \big] \right\} } .
\label{eq:new4}
\end{gather}
The second expression for $r_{1,2}$ is found from solving for the (1,2) entry of the residue of \eqref{eq:new3} at $x=b_6 t$
simultaneously with the (1,2) entry of the residue of \eqref{eq:new3} at $x=b_6^{-1}t$. This gives
\begin{gather}
 r_{1,2} = \frac{ - q t w (  b_6t  - \lambda) ( t - b_6 \lambda) }
 { t b_5 - t b_6 z_2 + b_5 b_6 ( b_6 t - \lambda) + t ( b_6 t - \lambda)	 (t - b_6 \lambda) }  .
\label{eq:new5}
\end{gather}
Combining \eqref{eq:new6} and \eqref{eq:new4} or \eqref{eq:new6} and \eqref{eq:new5}, and employing the change of variables
\eqref{Zsubs:a} and \eqref{Zsubs:b} with \eqref{FGsubs:a} and \eqref{FGsubs:b}, we can solve for $\hat{w}$ in two ways.
Assuming $w$ is non-zero it cancels out, leaving an expression for $\hat{\mathpzc{Z}}$ in terms of $\mathpzc{Z}$ and $\lambda$.
This is equivalent to the f\/irst $E^{(1)}_6$ $q$-Painlev\'e equation \eqref{1stEvol}.

To f\/ind the second equation we solve \eqref{eq:new3} for $\tilde{A}(x;t)$ $\vphantom{\hat{A}^{2^2}}$
\begin{gather*}
\tilde{A}(x;t) = \tilde{B}(qx;t)^{-1} \tilde{A}(x;qt) \tilde{B}(x;t),
%\label{eq:new7}
\end{gather*}
and use this to f\/ind the zero of $ \tilde{A}(x;t)_{12} $, i.e., $ g(t) $. In addition to $ r_{1,2} $
we now require $ r_{2,2} $ (even though the denominator of $ \tilde{A}(x;t)_{12} $ depends on $ r_{1,1}$, $r_{1,2}$, $r_{2,1}$, $r_{2,2} $
identities resulting from the compatibility conditions imply that this will trivialise~-- see the subsequent observation).
The entry $ r_{2,2} $ has already been found, along with $ r_{1,2}$, from the arguments given earlier and this is
\begin{gather}
 r_{2,2} =  \hat{\lambda}-\frac{qt\big(1+b_6^2\big)}{b_6}-\big(b_6qt-\hat{\lambda}\big)\big(qt-b_6\hat{\lambda}\big)\nonumber\\
 \hphantom{r_{2,2} = }{}
 \times \frac{b_6+q^2t^2 b_5\big(1+b_6^2\big)-b_5b_6qt\big(\hat{\alpha}+\hat{\lambda}\big)}
                                 {b_6^2 \big(\hat{\lambda}-b_6qt\big)-b_6qt+b_5b_6^2qt\hat{z}_1-b_5b_6qt\big(b_6qt-\hat{\lambda}\big)
                                 \big(qt-b_6\hat{\lambda}\big)} .
\label{r22}
\end{gather}
The numerator of $ \tilde{A}(x;t)_{12} $ appears to be a polynomial of degree 6 in $ x $, however it has trivial
zeros matching those of the denominator
\begin{gather*}
  \frac{q^2b_5}{b_6^2}(x-b_6t)(x-b_6qt)(b_6x-t)(b_6x-qt) \hat{w} ,
\end{gather*}
so that their ratio is in fact polynomial of degree~2.
Into $ \tilde{A}(x;t)_{12} $ we f\/irst
substitute for $ r_{1,2} $ using (\ref{eq:new6}), then for $ r_{2,2} $ using (\ref{r22}), and thirdly for
$ \hat{\alpha}$, $\hat{\beta}$, $\hat{\gamma}$, $\hat{\delta} $ using (\ref{alpha_new}), (\ref{beta_new}), (\ref{gamma_new}), (\ref{delta_new})
at up-shifted times, respectively. Into the resulting expression we employ (\ref{z_new}) for $ \hat{z}_1$, $\hat{z}_2 $
along with~(\ref{Zsubs:a}), (\ref{Zsubs:b}) at the up-shifted time to bring the whole expression in terms of
$ \hat{\lambda} $ and $ \hat{\mathpzc{Z}} $. The relevant zero of the ensuing expression (the other zero is $ x=0 $)
then gives $ \lambda=g $ in terms of $ \hat{\lambda}=\hat{g} $ and~$ f $, or equivalently by~(\ref{2ndEvol}).
\end{proof}

Now we recount the formulation given by Sakai \cite{Sa_2006} and Murata \cite{Mu_2009}. Their Lax pairs are
\begin{gather*}
Y(qx;t)  = A(x,t) Y(x;t), \qquad
Y(x;qt)  = B(x,t) Y(x;t),
\end{gather*}
satisfying the compatibility condition $A(x,qt)B(x,t) = B(qx,t)A(x,t)$.
The spectral matrix is parameterised in the following way
\begin{gather*}
A(x,t) = \begin{pmatrix} \kappa_1 W(x,t) & \kappa_2 w L(x,t) \\	
				    \kappa_1 w^{-1} X(x,t) & \kappa_2 Z(x,t) \end{pmatrix}
				     = A_0 + A_1x + A_2 x^2 + A_3 x^3,
\end{gather*}
subject to the key properties
\begin{enumerate}\itemsep=0pt
\item the determinant of $A(x,t)$ is
\[
\kappa_1 \kappa_2 (x-a_1)(x-a_2)(x-a_3)(x-a_4)(x-a_5t)(x-a_6t),
\]
\item $A_3$ is diagonal with entries $\kappa_1$ and $\kappa_2 = q\kappa_1$,
\item $A_0$ has eigenvalues $\theta_1 t$ and $\theta_2 t$,
\item the single root of the (1,2) entry of $A(x,t)$ in $x$ is $\lambda$,
\item $A(\lambda,t)$ is lower triangular with diagonal entries $\kappa_1 \mu_1$ and $\kappa_2 \mu_2$.
\end{enumerate}
Given these requirements, the entries of $A(x,t)$ are specif\/ied by
\begin{gather*}
L(x,t)  = x - \lambda, \\
Z(x,t)  =\mu _2+(x-\lambda ) \big[ \delta_2+x^2+x (\gamma +\lambda ) \big], \\
W(x,t)  = \mu _1+(x-\lambda ) \big[ \delta_1+x^2+x  (-\gamma -e_1 +\lambda  ) \big], \\
X(x,t)  = \big[ W Z - (x-a_1 ) (x-a_2 )  (x-a_3 )  (x-a_4 )  (x-a_5 t )
 (x-a_6  t  ) \big] L^{-1},
\end{gather*}
where
\begin{gather*}
(\kappa _1-\kappa _2) \delta_1  =  \lambda^{-1}  [ \kappa _1 \mu _1+\kappa _2 \mu _2 - \theta_1t -\theta_2 t  ]
	- \kappa _2 \big[ \gamma   (\gamma +e _1  ) + 2 \lambda ^2 - \lambda  e _1+ e _2 \big], \\
( \kappa _1-\kappa _2) \delta_2  =  - \lambda^{-1}  [ \kappa _1 \mu _1+\kappa _2 \mu _2-\theta_1t-\theta_2 t  ]
     + \kappa_1   \big[ \gamma  (\gamma + e_1 ) + 2 \lambda^2- \lambda  e_1+ e _2 \big] ,\\
\mu _1\mu_2  =   (\lambda -a_1 )  (\lambda -a_2 )  (\lambda -a_3 )  (\lambda -a_4 )
 (\lambda -a_5 t )  (\lambda -a_6 t ) ,\\
\theta_1 \theta_2  =   a_1 a_2 a_3 a_4 a_5 a_6q  \kappa _1^2 .
\end{gather*}
Here $e_j$ is the $j^{\rm th}$ elementary symmetric function of the indeterminates $\{ a_1,a_2,a_3,a_4,a_5t, a_6t \}$.
Despite the expression for $X(x,t)$, it is a quadratic polynomial in $x$. In Murata's notation we have $\tilde{\mu} = \mu_1$,
$\mu = \mu_2$, $\tilde{\delta} = \delta_1$ and $\delta = \delta_2$.
The deformation matrix $B(x,t)$ is a rational function in $x$ of the form
\begin{gather*}
B(x,t) = \frac{ x (x \mathbbmss{1} + B_0)} { (x-a_5 qt) ( x - a_6 qt)} .
\end{gather*}

Next we consider the f\/irst transformation of the Sakai linear problem with the following def\/inition:
\begin{gather*}
  {\mathpzc Y}(x,t) = \mathpzc{S}(x,t)^{-1} Y(x,t) ,
\end{gather*}
and
\begin{gather*}
 \mathpzc{S} = \begin{pmatrix}  1            & 0 \\	
			 s_1+ s_2x    & x
	       \end{pmatrix} .
\end{gather*}
The transformed spectral linear problem is
\begin{gather*}
{\mathpzc Y} (qx,t) = {\mathpzc A}(x,t) {\mathpzc Y}(x,t),
\end{gather*}
with a transformed spectral matrix
\begin{gather*}
{\mathpzc A}(x,t) = \mathpzc{S}(qx,t)^{-1} A(x,t) \mathpzc{S}(x,t) .
\end{gather*}
We f\/ix the parameters of the transformation by the requirement that the coef\/f\/icient of $x^{-1}$ in the (2,1) entry of $ \mathpzc{S} $
is zero (only the (2,1) entry is non-zero) and also that the coef\/f\/icient of $x^0$ in the (2,1) entry of $ \mathpzc{S} $ is zero. Thus
we f\/ind
\begin{gather*}
s_1 = \frac{1}{2q \kappa_1 w \lambda} \big[ (\theta_2-\theta_1)t + \kappa_1 ( \mu_1 - q \mu_2 + (q \delta_2 - \delta_1) \lambda  ) \big], \\
 q(q\theta_1 - \theta_2 ) w s_2  =
	\left(  \frac{2 \theta_1 \theta_2}{\kappa_1}t-q \theta_1 \mu_2 - \theta_2 \mu_1 \right)\frac{1}{\lambda^2}
	 	-q \kappa_1 e_5 t^{-1}\frac{1}{\lambda} \\
\hphantom{q(q\theta_1 - \theta_2 ) w s_2  =}{}
- e_1 \theta_2 +(q\theta_1 - \theta_2)\gamma + (q \theta_1 + \theta_2 ) \lambda .
\end{gather*}
The new spectral matrix can be parameterised by the polynomial
\begin{gather*}
{\mathpzc A}(x,t) = {\mathpzc A}_0 + {\mathpzc A}_1 x +  {\mathpzc A}_2 x^2 +  {\mathpzc A}_3 x^3,
\end{gather*}
and possesses the following properties
\begin{enumerate}\itemsep=0pt
\item the determinant of $\mathcal A(x,t)$ is
\[
\kappa_1^2 (x-a_1)(x-a_2)(x-a_3)(x-a_4)(x-a_5t)(x-a_6t),
\]
\item ${\mathpzc A}_3= \kappa_1 \mathbbmss{1}$,
\item ${\mathpzc A}_0$ is diagonal with entries $\theta_1 t$ and $q^{-1} \theta_2 t $,
\item the roots of the (1,2) entry of $\mathpzc A$ in $x$ are 0 and $\lambda$,
\item ${\mathpzc A}(\lambda,t)$ is lower triangular with diagonal entries $\kappa_1 \mu_1$ and $\kappa_2\mu_2$.
\end{enumerate}
Any such matrix admits the general form
\begin{gather*}
{\mathpzc A}(x,t) =
  \begin{pmatrix}
     \theta_1 t +  \kappa _1 x  [ (x-\lambda )  (x- \mathpzc{a}  )+\nu _1 ] &
     q \kappa _1 w x (x-\lambda )  \\
     \kappa _1 w^{-1} x  (x \mathpzc{c}+ \mathpzc{d}  )  &
     q^{-1}\theta_2 t + \kappa _1 x  [ (x-\lambda )  (x- \mathpzc{b}  )+\nu _2 ]
  \end{pmatrix} ,
\end{gather*}
where the properties given above f\/ix the introduced parameters as
\begin{gather*}
  (q\theta_1 - \theta_2) \mathpzc{a}   =
	  \big[ \theta_2 \nu _1+ q\theta_1 \nu _2 + q \kappa_1e_5t^{-1} \big]\lambda^{-1} + q\theta_1 e_1
		- 2q \theta_1 \lambda , \\
  (q\theta_1 - \theta_2) \mathpzc{b}   =
	- \big[ \theta_2 \nu _1+ q\theta_1 \nu _2 + q \kappa_1e_5t^{-1} \big]\lambda^{-1} - \theta_2 e_1
		+ 2 \theta_2 \lambda , \\ 		
  q\mathpzc{c}   =  \mathpzc{a}\mathpzc{b} + 2(\mathpzc{a}+\mathpzc{b})\lambda + \lambda^2 - e_2 + \nu_1 + \nu_2, \\
  q\mathpzc{d}   = - (\mathpzc{a}+\mathpzc{b})\lambda^2 - 2\mathpzc{a}\mathpzc{b} \lambda - \mathpzc{a}\nu_2 - \mathpzc{b}\nu_1
				+ ( q \mathpzc{c} - \nu_1 - \nu_2 )\lambda + e_3 + \frac{q\theta_1 + \theta_2}{q\kappa_1} t.
\end{gather*}
The variables, $\nu_1$ and $\nu_2$ are def\/ined by
\[
\nu_1 = \frac{\kappa_1\mu_1 - \theta_1 t }{\kappa_1\lambda}  \qquad \mbox{and} \qquad
	 \nu_2 = \frac{q\kappa_1\mu_2 - \theta_2 t}{q \kappa_1 \lambda}.
\]
The transformed deformation matrix $ \mathpzc{B} $ is computed using
\[
\mathpzc{B}(x,t) = \mathpzc{S}(x,qt)^{-1} B(x,t) \mathpzc{S}(x,t),
\]
and has the form
\[
\mathpzc{B} = \frac{x(x \mathpzc{B}_{0} + \mathbbmss{1}) } {(x- a_5 q t)( x - a_6 q t)} .
\]
We def\/ine a new variable $\nu$ using
\begin{gather*}
 \mu_2 \equiv \frac{(\lambda - a_1)(\lambda - a_2)(\lambda - a_3)(\lambda - a_4) }  { \lambda - \check{\nu}  } ,
\end{gather*}
and by implication
\begin{gather*}
 \mu_1 \equiv (\lambda - a_5 t) ( \lambda - a_6 t) ( \lambda - \check{\nu}) .
\end{gather*}
Using identical techniques to those employed in the proof of Lemma~\ref{lem:new1}, we can show that the compatibility
relation leads to the evolution equations
\begin{gather*}
(\lambda - \check{\nu})(\lambda - \nu)   = \dfrac{(\lambda-a_1)(\lambda-a_2)(\lambda-a_3)(\lambda-a_4)}
	{(\lambda - a_5 t)  (\lambda -  a_6t)},
%\label{eq:check1}
\\
\left(1-\dfrac{\nu}{\hat{\lambda}} \right)\left( 1- \dfrac{\nu}{\lambda}\right)  = \dfrac{a_5a_6}{q} \dfrac{(\nu-a_1)(\nu-a_2)(\nu-a_3)(\nu-a_4)}
                                                                                           {(a_5a_6t \nu + \theta_1/q\kappa_1)(a_5a_6t \nu + \theta_2/q\kappa_1)}.
%\label{eq:check2}
\end{gather*}

To make the full correspondence between our system and this one we must consider a further transformation,
given by the linear solution
\[
\mathfrak{Y}(x,t) = \frac{ \big[\vartheta_q \big( q^{-1} x  \big) \big]^3} { e_{q,t}(x)}  \mathpzc{Y}\big(x^{-1}, t^{-1}\big).
\]
The prefactors are elliptic functions def\/ined in terms of the $q$-factorial by
\begin{gather*}
\vartheta_q(z) = \big( q, -qz, -z^{-1}; q \big)_{\infty}, \qquad
    e_{q,t}(z) = \frac{ \vartheta_q (z) \vartheta_q\big(t^{-1}\big)} { \vartheta_q \big(z t^{-1} \big) },
\end{gather*}
with properties
\begin{gather*}
\vartheta_q (qz) = qz \vartheta_q(z), \qquad e_{q,t}(qz) = t e_{q,t}(z), \qquad e_{q,qt}(z) = z e_{q,t}(z).
\end{gather*}
This is the solution satisfying the linear equations
\begin{gather*}
\mathfrak{Y}\big(q^{-1} x,t\big) = \mathfrak{A}(x,t) \mathfrak{Y}(x,t) ,  \qquad
\mathfrak{Y}\big(x,q^{-1}t\big) = \mathfrak{B}(x,t) \mathfrak{Y}(x,t) .
\end{gather*}
The transformed spectral matrix $\mathfrak{A}$ is given by
\begin{gather}
\mathfrak{A}(x,t)  =  t x^3 \mathcal{A}\big(x^{-1}, t^{-1}\big) = \mathfrak{A}_3 + \mathfrak{A}_2 x + \mathfrak{A}_1 x^2 + \mathfrak{A}_0 x^3,
\label{eq:transSpec3}
\end{gather}
which swaps the roles of the leading matrices around $x = 0$ and $x= \infty$. This spectral matrix has the properties
\begin{enumerate}\itemsep=0pt
\item the determinant of $\mathfrak A(x,t)$ is
\[
\kappa_1^2 (1-a_1x)(1-a_2x)(1-a_3x)(1-a_4x)(t-a_5x)(t-a_6x),
\]
\item $\mathfrak A_3 =  \kappa_1 t \mathbbmss{1} $,
\item $\mathfrak{A}_0$ is diagonal with entries $\theta_1$ and $q^{-1} \theta_2$,
\item the roots of the (1,2) entry of $\mathfrak A(x,t)$ in $x$ are 0 and $\lambda^{-1}$,
\item $\mathfrak A(\lambda,t)$ is lower triangular with diagonal entries $ \kappa_1 \mu_1 t \lambda^{-3}$ and
$ \kappa_1 \mu_2 t \lambda^{-3}$.
\end{enumerate}
The transformed deformation matrix has the form
\begin{gather}
\mathfrak{B}(x,t) = \frac{x} {(t - a_5 qx) (t - a_6 qx)} (x\mathbbmss{1} + \mathfrak{B}_{0}) .
\label{eq:transDef4}
\end{gather}
Since the compatibility relation between \eqref{eq:transSpec3} and \eqref{eq:transDef4} is rationally equivalent to that for
$ \mathpzc{Y} $, the evolution equations are the same.

It is clear that $\mathfrak{Y}$ and $Y_n$ satisfy equivalent linear problems and that the following
correspondences hold:
\begin{gather*}
q   \mapsto q^{-1} ,   \qquad
t    \mapsto t^{-1} ,  \qquad
\lambda(t)   \mapsto \frac{1}{g(t^{-1})} , \qquad
\nu(t)   \mapsto f\big(t^{-1}\big) ,
\\
\kappa_1 \mapsto b_6, \qquad   a_i \mapsto b_i  \quad i = 1,2,3,4 , \qquad a_5 \mapsto \frac{1}{b_6} , \qquad
\theta_1 \mapsto - b_5 b_6, \qquad  q^{-1} \theta_2  \mapsto - \frac{b_6}{b_5} .
\end{gather*}

\subsection{Reconciliation with the Lax pair of Yamada \cite{Ya_2011}}

In his derivation of a Lax pair for the $E_6^{(1)}$ $q$-Painlev\'e system Yamada employed the degeneration limits
of $E_8^{(1)}$ $q$-Painlev\'e $\rightarrow$ $E_7^{(1)}$ $q$-Painlev\'e $\rightarrow$ $E_6^{(1)}$ $q$-Painlev\'e.
In doing so he retained eight parameters $b_1,\ldots,b_8$ constrained by $q b_1 b_2 b_3 b_4 = b_5 b_6 b_7 b_8$, and his
$E_6^{(1)}$ $q$-Painlev\'e equation was given by the mapping of the variables
\begin{gather*}
t \mapsto q^{-1} t, \qquad f,g \mapsto \bar{f}, \bar{g},
\end{gather*}
subject to the coupled f\/irst-order system (see his (36))
\begin{gather}
\frac{(fg-1) (\bar{f} g -1) } {f \bar{f} }  = q \frac{ (b_1 g -1) (b_2g -1) (b_3g -1) (b_4g -1) } { b_5 b_6 (b_7 g -t) (b_8 g -t)},
\label{Y_E6qP:a}  \\
\frac{ (f \underline{g} -1) (fg -1) } {g \underline{g} }  = \frac{ (b_1 -f) (b_2 -f ) (b_3 -f) (b_4 -f) } { (f - b_5t) (f - b_6 t) }.
\label{Y_E6qP:b}
\end{gather}
The Lax pairs constructed by the degeneration limits were given as a coupled second-order $q$-dif\/ference equation in a scalar
variable $Y(z,t)$ (see his~(37))
\begin{gather}
\frac{(b_1 q  -z) (b_2 q -z) (b_3 q -z) (b_4 q -z) t^2  } {q(qf -z) z^4 } \left[ Y(q^{-1}z) - \frac{g z}{t^2 (gz- q)} Y(z) \right] \nonumber\\
\qquad{} + \left[ \frac{q(b_1g -1) (b_2 g -1) (b_3g -1) (b_4 g -1) } {g (fg -1) z^2 (gz -q)}
		- \frac{b_5 b_6 (b_7g -t) (b_8 g -t) } { f g z^3  } \right] Y(z) \nonumber\\
\qquad{} + \frac{ (b_5t -z) (b_6 t -z)   }   { t^2 z^2 (f-z)  } \left[ Y(qz) - \frac{t^2 (gz -1) } {gz} Y(z) \right] = 0,
\label{Y_E6LP:a}
\end{gather}
and a second-order, mixed $q$-dif\/ference equation,
\begin{gather}
\frac{g z} {t^2} Y(z) + (q -g z) Y\big(q^{-1} z\big) - q^{-2} gz(qf -z) \bar{Y}\big(q^{-1} z\big) = 0.
\label{Y_E6LP:b}
\end{gather}

In order to bring \eqref{Y_E6qP:a} and \eqref{Y_E6qP:b} into correspondence with our form of the $E_6^{(1)}$ $q$-Painlev\'e
system (see \eqref{eq:intro1} and \eqref{eq:intro2}) we will employ the following transformation of Yamada's
variables
\begin{gather*}
t \mapsto t^{-1}, \qquad z \mapsto z^{-1}, \qquad
f,g \mapsto g^{-1}, f^{-1}, \qquad\tilde{Y}(z) = Y\big(z^{-1}\big),
\end{gather*}
and the specialisations of the parameters
\begin{gather*}
b_5 \mapsto b_6^{-1}, \qquad b_6 \mapsto b_6 , \qquad
b_7 \mapsto q b_5, \qquad b_8 \mapsto b_5^{-1},
\end{gather*}
so that $b_5 b_6 =1$ and $b_7 b_8 =q$. Under these transformations we deduce that~\eqref{Y_E6qP:a} becomes~\eqref{eq:intro2}
and~\eqref{Y_E6qP:b} becomes~\eqref{eq:intro1}. Furthermore the pure second-order divided-dif\/ference equation~\eqref{Y_E6LP:a}
becomes
\begin{gather}
\frac{ \prod\limits_{j=1}^4 (1-b_j z)}{t^2z(z - g)} \tilde{Y}(z)
+ \left\{ - \frac{ \prod\limits_{j=1}^4 (1-b_j z)}{z(z - g)(1-f z)} + \frac{z \prod\limits_{j=1}^4 (b_j - f )}{f (1 - fg) (1-f z)}
          - \frac{ z (f -b_5 q t) ( b_5 f - t)} {b_5 q t^2 f } \right.\nonumber\\[-1ex]
\left. \qquad{} - \frac { (z -b_6 q t) (b_6z -qt) (q -fz) } { b_6qt^2 z (z -qg) }\vphantom{\frac{ \prod\limits_{j=1}^4 (1-b_j z)}{z(z - g)(1-f z)}} \right\}
		\tilde{Y}\big(q^{-1}z\big)
			+ \frac{ (z - b_6qt) (b_6z -qt)} {b_6 z (z-qg) } \tilde{Y}\big(q^{-2}z\big) =0,\!\!\!
\label{eq:NEW2ndODD}
\end{gather}
and the mixed divided-dif\/ference equation \eqref{Y_E6LP:b} becomes
\begin{gather}
\frac{qt^2}{fz} \tilde{Y}\big(q^{-1}z;t\big) - q \frac{(1-fz)}{fz} \tilde{Y}(z;t) - \frac{(z-g)} {f g z^2} \tilde{Y}(z;qt) = 0 .
\label{eq:NEWdivdiff}
\end{gather}

Having put Yamada's Lax pairs into a suitable form we now seek to make a correspondence with our own theory and results.
A single mixed divided-dif\/ference equation can be constructed from the matrix Lax pairs (\eqref{spectral_DDO} and \eqref{deform_DDO}). For generic semi-classical
systems on a $q$-lattice grid we can deduce either
\begin{gather*}
- \frac{1} {W + \Delta y V} \frac{1} {\mathfrak P_+} p_n(x;qt) + \frac{1} {(W + \Delta y V) ( R + \Delta u S) }
	\left[ - \frac{\mathfrak W_+} {\mathfrak T_+} + \frac { \mathfrak R_+ } { \mathfrak P_+ } \right] p_n(x;t)\\
\qquad{}
	+ \frac{1} {R + \Delta u S } \frac{1} {\mathfrak T_+} p_n(qx;t) =0,
\end{gather*}
or an alternative,
\begin{gather}
 - \frac{\mathfrak T_+(x)}{(W - \Delta y V)(x)} p_n(qx;qt)
	+ \frac{1}{(W - \Delta y V)(x)(R + \Delta u S)(qx)}\nonumber\\
\qquad{}\times \big[ \mathfrak T_+(x) \mathfrak R _+(qx) + \mathfrak P_+(qx) \mathfrak W_-(x)  \big] p_n(qx;t)
		- \frac{\mathfrak P_+(qx)}{(R + \Delta u S)(qx)} p_n(x;t) =0,
\label{mixedDDE}
\end{gather}
which we will work with. Using the spectral and deformation data (\ref{eq:weight3}), (\ref{Ddata1}) and the explicit evaluations
of the deformation matrix (\ref{DMatrix11}) and spectral matrix (\ref{eq:Anstar2}), we compute the
coef\/f\/icients of the above equation
\begin{gather*}
 -(R + \Delta u S)(qx)\mathfrak T_+(x) = a_n\frac{1-qb_5^2}{b_5}x(x-b_6t)(x-g) ,
\\
 -(W - \Delta y V)(x)\mathfrak P_+(qx) = \frac{a_n\gamma_n}{\hat{\gamma}_n}\frac{b_6\big(1-qb_5^2\big)}{b_5}\frac{t(1-b_2x)(1-b_3x)(x-b_6t)}{b_5qt-f} ,
\\
  \mathfrak T_+(x) \mathfrak R _+(qx) + \mathfrak P_+(qx) \mathfrak W_-(x)
  = \frac{a_n\gamma_n}{\hat{\gamma}_n}\frac{b_6\big(1-qb_5^2\big)}{b_5}\frac{(x-b_6t)(b_6x-t)(1-fx)}{b_5qt-f} .
\end{gather*}
Now we set $ p_n = FU $ where $ F $ is a gauge factor and $ U $ is the new independent variable, into~(\ref{mixedDDE}) and make a direct comparison with~(\ref{eq:NEWdivdiff}). Comparing the coef\/f\/icients of
$ U(x;t) $ and $ U(qx;t) $ in this later equation we deduce that
\begin{gather*}
  \frac{F(qx,t)}{F(x,t)} = \frac{1}{t^2}\frac{(1-b_2x)(1-b_3x)}{1-b_6xt^{-1}} .
\end{gather*}
A solution is given by
\begin{gather*}
  F(x,t) = e_{q,t^{-2}}(x)\frac{\big(b_6xt^{-1};q\big)_{\infty}}{(b_2x,b_3x;q)_{\infty}} C(x,t) ,
\end{gather*}
where $ C $ is a $q$-constant function, $ C(qx,t)=C(x,t) $. Now comparing the coef\/f\/icients of $ U(qx;qt) $ and
$ U(qx;t) $ in the previous equation we f\/ind that
\begin{gather*}
  \frac{F(qx,qt)}{F(qx,t)} = \frac{\gamma_n}{\hat{\gamma}_n}\frac{b_6(b_6x-t)}{qg(b_5qt-f)x^2} .
\end{gather*}
Substituting our solution into this equation we f\/ind a complete cancellation of all the $x$
dependent factors resulting in a pure $q$-dif\/ference equation in $ t $
\begin{gather*}
 \frac{\hat{\gamma}_n\hat{C}}{\gamma_nC} = \frac{b_6qt}{g(f-b_5qt)} .
\end{gather*}
Thus we just need a solution $ C(t) $ independent of $x$, however we only require the existence of
a~non-zero, bounded solution rather than knowledge of a specif\/ic solution. In conclusion we f\/ind that
our new mixed, divided-dif\/ference equation is now{\samepage
\begin{gather*}
 t^2 U(x;t) - (1-fx)U(qx;t) - \frac{x-g}{qgx}U(qx;qt) = 0 ,
\end{gather*}
which is clearly proportional to (\ref{eq:NEWdivdiff}) with the identif\/ication $ U(x;t)=\tilde{Y}(q^{-1}x;t) $.}

A second-order $q$-dif\/ference equation in the spectral variable $x$ for one of the components, say $p_n$,
was given in \eqref{eq:LSODDE}, and for $q$-linear grids can be simplif\/ied as
\begin{gather}
\frac{W + \Delta y V}{\mathfrak T_+}(x) p_n(q x) - \left[ \frac{ \mathfrak{W}_+ } {\mathfrak T_+} (x) + \frac{ \mathfrak{W}_{-} } {\mathfrak T_+} \big(q^{-1}x\big) \right] p_n(x)\nonumber\\
\qquad{}
	+ \frac{W - \Delta y V}{\mathfrak T_+} \big(q^{-1}x\big) p_n\big(q^{-1}x\big) = 0.
\label{q-LSODDE}
\end{gather}
From the explicit solution of the gauge factor we note
\begin{gather*}
  \frac{F(qx,t)}{F(x,t)} = \frac{(1-b_2x)(1-b_3x)}{t(t-b_6x)}, \qquad
  \frac{F\big(q^{-1}x,t\big)}{F(x,t)} = \frac{t\big(t-b_6q^{-1}x\big)}{\big(1-b_2q^{-1}x\big)\big(1-b_3q^{-1}x\big)} .
\end{gather*}
Substituting the change of variables into (\ref{q-LSODDE}) we compute that
\begin{gather*}
 \frac{W + \Delta y V}{\mathfrak T_+}(x)\frac{F(qx,t)}{F(x,t)} =
 \frac{1}{(q-1)u_1 a_n}\frac{b_6}{t}\frac{\prod\limits^{4}_{j=1}(1-b_jx)}{x(x-g)} ,
\\
 \frac{W - \Delta y V}{\mathfrak T_+}\big(q^{-1}x\big)\frac{F\big(q^{-1}x,t\big)}{F(x,t)} =
 \frac{1}{(q-1)u_1 a_n}\frac{t(x-b_6 qt)(b_6 x-qt)}{x(x-qg)} .
\end{gather*}
In addition, using the explicit representations of the diagonal elements of $A_n^*$, i.e., $\mathfrak W_{\pm}$
(see \eqref{eq:Anstar1}, \eqref{eq:Anstar2}) we compute that
\begin{gather*}
 - \frac{1}{b_6 t} \left[ \frac{\mathfrak W_{+} } { x(x-g) } + \left. \frac{ \mathfrak W_- } { x(x-g) } \right|_{q^{-1}x} \right]
  = - \frac{ \prod\limits_{j=1}^4 (1 -b_j x) } {x(x-g) (1-fx) } + \frac{x \prod\limits_{j=1}^4 (f -b_j) } {f(1-fg) (1-fx) }\nonumber\\
 \qquad{}
    - \frac{x(f -b_5qt) (b_5f -t) } { b_5 q t^2 f} - \frac{(x -b_6qt) (b_6x -qt) (q-fx) }{b_6qt^2 x(x-qg) }.
%\label{eq:Z2}
\end{gather*}
In summary we f\/ind
\begin{gather}
 \frac{ \prod\limits_{j=1}^4 (1-b_j x) } {t^2 x(x-g) } U(qx)
  + \left[
- \frac{ \prod\limits_{j=1}^4 (1 -b_j x) } {x(x-g) (1-fx) } + \frac{x \prod\limits_{j=1}^4 (f -b_j) } {f(1-fg) (1-fx) } \right.
 - \frac{x(f -b_5qt) (b_5f -t) } { b_5 q t^2 f}\nonumber\\
    \left.
    \qquad{} - \frac{(x -b_6qt) (b_6x -qt) (q-fx) }{b_6qt^2 x(x-qg) }\vphantom{\frac{ \prod\limits_{j=1}^4 (1 -b_j x) } {x(x-g) (1-fx) }}
    \right] U(x)
		+ \frac{(x-b_6qt) (b_6x -qt) } {b_6 x(x-qg) } U\big(q^{-1} x\big) = 0 .		
\label{eq:Final2ndOdivdiff}
\end{gather}
Thus we can see that \eqref{eq:Final2ndOdivdiff} agrees exactly with \eqref{eq:NEW2ndODD} and the
identif\/ication noted above.

\subsection*{Acknowledgements}
This research has been supported by the Australian Research Council's Centre of Excellence for Mathematics
and Statistics of Complex Systems. We are grateful for the clarif\/ications by Kenji Kajiwara of results given
in~\cite{KMNOY_2004} and~\cite{KMNOY_2005} and the assistance of Yasuhiko Yamada in explaining the results
of his work~\cite{Ya_2011}. We also appreciate the assistance of Jason Whyte in the preparation of this manuscript.

\pdfbookmark[1]{References}{ref}

\LastPageEnding


\begin{thebibliography}{99}
\footnotesize\itemsep=0pt

\bibitem{Birkhoff_1911}
Birkhof\/f G.D., General theory of linear dif\/ference equations, \href{http://dx.doi.org/10.2307/1988577}{\textit{Trans.
  Amer. Math. Soc.}} \textbf{12} (1911), 243--284.

\bibitem{Birkhoff_1913}
Birkhof\/f G.D., The generalized {R}iemann problem for linear dif\/ferential
  equations and the {A}llied problems for linear dif\/ference and $q$-dif\/ference
  equations, \href{http://www.jstor.org/stable/20025482}{\textit{Trans. Amer. Math. Soc.}} \textbf{49} (1913), 521--568.

\bibitem{Ismail_2005}
Ismail M.E.H., Classical and quantum orthogonal polynomials in one variable,
  \textit{Encyclopedia of Mathematics and its Applications}, Vol.~98, Cambridge
  University Press, Cambridge, 2005.

\bibitem{JS_1996}
Jimbo M., Sakai H., A {$q$}-analog of the sixth {P}ainlev\'e equation,
  \href{http://dx.doi.org/10.1007/BF00398316}{\textit{Lett. Math. Phys.}} \textbf{38} (1996), 145--154,
  \href{http://arxiv.org/abs/chao-dyn/9507010}{chao-dyn/9507010}.

\bibitem{KMNOY_2005}
Kajiwara K., Masuda T., Noumi M., Ohta Y., Yamada Y., Construction of
  hypergeometric solutions to the {$q$}-{P}ainlev\'e equations, \href{http://dx.doi.org/10.1155/IMRN.2005.1439}{\textit{Int.
  Math. Res. Not.}} \textbf{2005} (2005), no.~24, 1439--1463,
  \href{http://arxiv.org/abs/nlin.SI/0501051}{nlin.SI/0501051}.


\bibitem{KMNOY_2004}
Kajiwara K., Masuda T., Noumi M., Ohta Y., Yamada Y., Hypergeometric solutions
  to the {$q$}-{P}ainlev\'e equations, \href{http://dx.doi.org/10.1155/S1073792804140919}{\textit{Int. Math. Res. Not.}}
  \textbf{2004} (2004), no.~47, 2497--2521, \href{http://arxiv.org/abs/nlin.SI/0403036}{arXiv:nlin.SI/0403036}.


\bibitem{KLS_2010}
Koekoek R., Lesky P.A., Swarttouw R.F., Hypergeometric orthogonal polynomials
  and their {$q$}-analogues, \href{http://dx.doi.org/10.1007/978-3-642-05014-5}{\textit{Springer Monographs in Mathematics}},
  Springer-Verlag, Berlin, 2010.

\bibitem{KS_1998}
Koekoek R., Swarttouw R.F., The Askey-scheme of hypergeometric orthogonal
  polynomials and its $q$-analogue, Report 98-17, Faculty of Technical
  Mathematics and Informatics, Delft University of Technology, 1998,
  \url{http://aw.twi.tudelft.nl/~koekoek/askey/}.

\bibitem{Ma_1988}
Magnus A.P., Associated {A}skey--{W}ilson polynomials as {L}aguerre--{H}ahn
  orthogonal polynomials, in Orthogonal Polynomials and their Applications
  ({S}egovia, 1986), \href{http://dx.doi.org/10.1007/BFb0083366}{\textit{Lecture Notes in Math.}}, Vol.~1329, Springer,
  Berlin, 1988, 261--278.

\bibitem{Ma_1995}
Magnus A.P., Special nonuniform lattice (snul) orthogonal polynomials on
  discrete dense sets of points, \href{http://dx.doi.org/10.1016/0377-0427(95)00114-X}{\textit{J.~Comput. Appl. Math.}} \textbf{65}
  (1995), 253--265, \href{http://arxiv.org/abs/math.CA/9502228}{math.CA/9502228}.

\bibitem{Mu_2009}
Murata M., Lax forms of the {$q$}-{P}ainlev\'e equations, \href{http://dx.doi.org/10.1088/1751-8113/42/11/115201}{\textit{J.~Phys.~A:
  Math. Theor.}} \textbf{42} (2009), 115201, 17~pages, \href{http://arxiv.org/abs/0810.0058}{arXiv:0810.0058}.

\bibitem{PNGR_1992}
Papageorgiou V.G., Nijhof\/f F.W., Grammaticos B., Ramani A., Isomonodromic
  deformation problems for discrete analogues of {P}ainlev\'e equations,
  \href{http://dx.doi.org/10.1016/0375-9601(92)90905-2}{\textit{Phys. Lett.~A}} \textbf{164} (1992), 57--64.

\bibitem{RGTT_2001}
Ramani A., Grammaticos B., Tamizhmani T., Tamizhmani K.M., Special function
  solutions of the discrete {P}ainlev\'e equations, \href{http://dx.doi.org/10.1016/S0898-1221(01)00180-8}{\textit{Comput. Math.
  Appl.}} \textbf{42} (2001), 603--614.


\bibitem{Sa_2005}
Sakai H., A {$q$}-analog of the {G}arnier system, \href{http://dx.doi.org/10.1619/fesi.48.273}{\textit{Funkcial. Ekvac.}}
  \textbf{48} (2005), 273--297.

\bibitem{Sa_2006}
Sakai H., Lax form of the {$q$}-{P}ainlev\'e equation associated with the
  {$A\sp {(1)}\sb 2$} surface, \href{http://dx.doi.org/10.1088/0305-4470/39/39/S13}{\textit{J.~Phys.~A: Math. Gen.}} \textbf{39}
  (2006), 12203--12210.

\bibitem{Sa_2001}
Sakai H., Rational surfaces associated with af\/f\/ine root systems and geometry of
  the {P}ainlev\'e equations, \href{http://dx.doi.org/10.1007/s002200100446}{\textit{Comm. Math. Phys.}} \textbf{220} (2001),
  165--229.

\bibitem{Wi_2010a}
Witte N.S., Semi-classical orthogonal polynomial systems on non-uniform
  lattices, deformations of the {A}skey table and analogs of isomonodromy,
  \href{http://arxiv.org/abs/1204.2328}{arXiv:1204.2328}.

\bibitem{Ya_2011}
Yamada Y., Lax formalism for {$q$}-{P}ainlev\'e equations with af\/f\/ine {W}eyl
  group symmetry of type {$E^{(1)}_n$}, \href{http://dx.doi.org/10.1093/imrn/rnq232}{\textit{Int. Math. Res. Not.}}
  \textbf{2011} (2011), no.~17, 3823--3838, \href{http://arxiv.org/abs/1004.1687}{arXiv:1004.1687}.

\end{thebibliography}
\end{document}